\documentclass[11pt,a4paper,reqno]{amsart} 
\usepackage{tikz}
\usetikzlibrary{positioning,decorations.pathmorphing}

\usepackage{amssymb,graphicx}

\usepackage{amssymb}

\newcommand{\koniec}{\begin{flushright}  $\Box $ \end{flushright}}
\newtheorem{theo}{Theorem}[section] 
\newtheorem{prop}[theo]{Proposition}  
\newtheorem{lemma}[theo]{Lemma}

\def\theequation{\thesection.\arabic{equation}}

\newcounter{mnotecount}[section]

\renewcommand{\themnotecount}{\thesection.\arabic{mnotecount}}

\newcommand{\mnote}[1]
{\protect{\stepcounter{mnotecount}}$^{\mbox{\footnotesize
$
\bullet$\themnotecount}}$ \marginpar{
\raggedright\tiny\em
$\!\!\!\!\!\!\,\bullet$\themnotecount: #1} }
\newcommand{\md}[1]{\mnote{{\bf MD:}#1}}
\newcommand{\hook}{{\setlength{\unitlength}{11pt}   
                   \begin{picture}(.833,.8)
                   \put(.15,.08){\line(1,0){.35}}
                   \put(.5,.08){\line(0,1){.5}}
                   \end{picture}}}
\newcommand{\CP}{\mathbb{CP}}
\newcommand{\C}{\mathbb{C}}

\newcommand{\Z}{\mathbb{Z}}

\newcommand{\R}{\mathbb{R}}

\def\p{\partial}

\def\be{\begin{equation}}

\def\ee{\end{equation}}

\def\bea{\begin{eqnarray}}
\def\eea{\end{eqnarray}}

\newcommand{\spp}{\mathbb{S}}

\begin{document}\date{13 May  2021}
\vspace*{-1.0cm}
\title{Conformal geodesics on gravitational instantons}
\author{Maciej Dunajski}
\address{Department of Applied Mathematics and Theoretical Physics\\ 
University of Cambridge\\ Wilberforce Road, Cambridge CB3 0WA, UK.}
\email{m.dunajski@damtp.cam.ac.uk}

\author{Paul Tod}
\address{The Mathematical Institute\\
Oxford University\\
Woodstock Road, Oxford OX2 6GG\\ UK.
}
\email{tod@maths.ox.ac.uk}
\begin{abstract} 
We study the integrability of the conformal geodesic flow
(also known as the conformal circle flow)
on the $SO(3)$--invariant gravitational instantons.  On a hyper--K\"ahler four--manifold the conformal 
geodesic equations
reduce to geodesic equations of a charged particle moving in a constant
self--dual magnetic field. In the case of the anti--self--dual Taub NUT
instanton we integrate these equations completely by separating 
the Hamilton--Jacobi equations, and finding a commuting set of first integrals. 
This gives the first example of an integrable conformal geodesic flow on a four--manifold which is not a symmetric space.
In the case of the 
Eguchi--Hanson we find all conformal geodesics
which lie on the three--dimensional 
orbits of the isometry group. In the non--hyper--K\"ahler case of 
the Fubini--Study metric on $\CP^2$
we use the first integrals arising from the conformal Killing--Yano tensors to recover the known complete integrability of conformal geodesics.
\end{abstract}   
\maketitle
\section{Introduction}
The geodesic flow on a (pseudo) Riemannian manifold $(M, g)$ is integrable if the underlying metric admits a sufficient number
of Killing vectors, or Killing tensors. The former correspond to first integrals linear in the velocity, and the latter
give polynomial first integrals of higher order. For example the geodesic motion on a round two--sphere is integrable as
any generator of the isometry group $SO(3)$ gives a Killing vector commuting with the Hamiltonian. On a tri--axial
ellipsoid there are no Killing vectors, but the additional quadratic 
first integral required for integrability is given by a rank--two Killing
tensor \cite{ellipsoid}.

Let $u$ be a unit tangent vector to a curve $\Gamma\subset M$, and let
$a\equiv \nabla_u u$ be the acceleration of $\Gamma$. 
The geodesic condition $a=0$ is not invariant under 
conformal rescalings of the metric $g\rightarrow \Omega^2 g$, but there is a different preferred set of curves on manifolds endowed with only conformal structures $(M, [g])$, where $[g]=\{\Omega^2 g, \Omega:M\rightarrow \R^+$\}.
These {\em conformal geodesics} are also known as {\em conformal circles}, and arise as solutions to a system
of third order ODEs on $M$. A conformal geodesic $\Gamma$ is uniquely specified by a point, a tangent direction, and a perpendicular acceleration. If $g$ is a representative metric in the conformal class, and $\nabla$ is the Levi--Civita connection of $g$, then the conformal geodesic equations are
\be
\label{conf_circ}
\nabla_u a= -(|a|^2  +L(u, u))u+L^{\sharp}(u), \quad
\ee
where $L\in\Gamma(TM\otimes TM)$ is the Schouten tensor given in terms of the Ricci tensor $R$ and the Ricci scalar $S$ by
\be
\label{schouten_tn}
L=\frac{1}{n-2}\Big(R-\frac{1}{2(n-1)}Sg\Big),
\ee
and
$L^{\sharp}:TM\rightarrow TM$ is the endomorphism defined by $g(L^{\sharp}(X), Y)=L(X, Y)$ for all vector fields $X, Y$.
It can be demonstrated \cite{BE, tractor, tod_circles, sihlan}  that the conformal geodesics only depend on the conformal class
of $[g]$, and not on the choice of the representative metric.
In (\ref{schouten_tn}) it is assumed that $n>2$. The case where $n=2$ will be discussed in \S\ref{cccm}.

Neither the Killing vectors nor the conformal Killing vectors of $g$ give rise to first integrals of (\ref{conf_circ}), and
it is  natural to ask whether there are examples of integrable conformal geodesic motions,  and what geometric
structures on $(M, [g])$ give rise to this integrablity. 
In \cite{tod_circles, rodar}
it was shown that the conformal Killing--Yano two--forms (CKY) 
give rise to first integrals of (\ref{conf_circ}):  $Y\in \Lambda^2(T^*M)$ 
is a CKY
if ${(\nabla Y)}_0\in \Gamma (\Lambda^3(T^*M))$, where $T_0$ denotes the 
trace--free part of $T$. The corresponding first integral is then
\be
\label{first_int_int}
Q=Y(u, a)-\frac{1}{n-1} \mbox{div}(Y)(u),
\ee
where $n=\mbox{dim}(M)$, and $\mbox{div}=*d*$ is the divergence.
This was sufficient \cite{tod_circles} to integrate (\ref{conf_circ}) on a non--conformally flat ${\bf Nil}$ 3--manifold, and a squashed 3--sphere 
 but already the conformal class of the Schwarzchild metric proved too difficult to handle, as there do not exist sufficiently many CKYs. 

 In this paper we shall
study the integrability of (\ref{conf_circ}) on four--dimensional conformal structures corresponding to some gravitational instantons:
solutions to Einstein equations on Riemannian four--manifolds which are compact or complete metrics which asymptotically approach a locally flat space - the decay rate, as well as the topology at infinity varries between different examples of instantons. See \cite{GHr, Dbook} for details. 
 
We shall explore the existence
of three CKYs corresponding to the underlying hyper--K\"ahler structure, and reduce (\ref{conf_circ}) to a system
of 2nd order ODEs corresponding to a forced geodesic equation in a constant magnetic field. Further progress can be  made for the anti--self--dual (ASD) Taub--NUT and the Eguchi--Hanson instantons. In the ASD Taub--NUT case we shall establish complete integrability by separating the Hamilton--Jacobi equation.
\begin{theo}
\label{thm_TN}
The conformal geodesic equation on the anti--self--dual Taub--NUT manifold is completely integrable in the Arnold--Liouville sense:
it reduces to a geodesic equation in a self--dual magnetic field which admits four first integrals in involution. 

The associated Hamilton--Jacobi equation is separable, and the conformal geodesic equations reduce to quadratures.
\end{theo}
 The Eguchi--Hanson metric
admits an isometric and tri--holomorphic action
of $SO(3)$  which preserves the magnetic field, and therefore gives rise to 
three charged linear first integrals of the Lorentz force equation. While full integrablity cannot be established in this case, there is enough structue
to reduce the conformal geodesic motion to quadratures under an additional assumption that the conformal geodesics lie on the 3--dimensional orbits of the $SO(3)$ subgroup of the isometry group.

The paper is organised as follows. In the next section we shall introduce the conformal geodesic
equations (\ref{conf_circ}), and focus on the special case where the underlying conformal structure admits an Einstein metric.
We shall discuss the first--integrals of (\ref{conf_circ}), and establish a link with the Lorentz force equations for
hyper--K\"ahler four manifolds (Proposition \ref{prop31}).
In \S\ref{sectionTN}  and \S\ref{sectionEH}  we shall study  (\ref{conf_circ}) 
 on the ASD Taub--NUT metric, and the Eguchi--Hanson metric respectively.
 The ASD Taub--NUT case is completely integrable (Theorem \ref{thm_TN} will follow
from Proposition \ref{prop_TN1} and Proposition \ref{prop_TN2}),
and  the Eguchi--Hanson case has a couple of integrable sub--cases (Proposition
\ref{propEH}).
 
 In \S\ref{section3} we shall establish the integrability of
the conformal geodesic flow for 
the Fubini--Study metric on $\CP^2$. 
We shall establish it by making use of nine CKYs on this space, and recover the results of \cite{jap2}, where all conformal geodesics of $\CP^n$ have been characterised by their horizontal lifts to helices on the  total
space of the fibration $S^{2n+1}\rightarrow \CP^n$. 
In Appendix A we shall rule out the existence of non--flat Riemannian Gibbons--Hawking metrics with three commuting vector fields, and in Appendix B
we shall provide necessary and sufficient conditions for a Killing trajectory
to be a conformal geodesic. 
\subsection*{Acknowledgements.} The work of MD has been partially 
supported by STFC consolidated grants ST/P000681/1, and  ST/T000694/1.
\section{Conformal geodesics on Einstein manifolds}
\label{section2}
Let $(M, g)$ be a Riemannian manifold. When using the index-free notation
we shall denote the $g$--inner product of two vector fields $X$ and $Y$ by
$g(X, Y)$. We also set $|X|^2\equiv  g(X, X)$, and use notation
$X\hook\psi$ for the $(p-1)$--form arising as a contraction of the $p$--form
$\psi$ with the vector field $X$. If $\Gamma$ is a curve, and $u$ is
a tangent vector to $\Gamma$, then $\nabla_u$ denotes the directional derivative along $\Gamma$, where $\nabla$ is the Levi--Civita connection of $g$.

In explicit computations involving conformal geodesics, and conformal Killing--Yano tensors  it is convenient to adopt the abstract index notation  
\cite{PR}. 
Thus $u^a, a=1, 2, \dots, \mbox{dim}(M)$ denotes a vector, $u^a\nabla_a v^b$ denotes the directional covariant derivative $\nabla_u v$ of another vector $v$.
On a Riemannian manifold $(M, g)$ an
isomorphism between $TM$ and $T^*M$ is realised by $u_a=g_{ab}u^b$, where
now $u_a$ denotes a one--form, and the Einstein summation convention is used.
Despite the presence of these indices, no choice of basis has been made. 
\vskip5pt

If the metric $g$ is Einstein, i.e.
$R=(S/n)g $ then the conformal geodesic equation
(\ref{conf_circ}) for a curve $\Gamma$ parametrised by an
arc--lengh $s$, and with a unit tangent vector $u$
reduces
to
\be
\label{conf_circ_e}
\nabla_u a=-|a|^2 u, \quad \mbox{where}\quad
|a|^2\equiv g(a, a)\geq 0\,\,\mbox{is a constant}, \quad\mbox{and}\quad
g(a, u)=0.
\ee
This system of 3rd order equations (\ref{conf_circ_e}) has long been studied in 
Riemannian geometry
\cite{jap1, jap2}, where the corresponding solution curves have been called circles. The terminology is motivated by the observation
that a development\footnote{In references \cite{jap1, jap2} the  equation
(\ref{conf_circ_e}) is written as
$
\nabla_s X_s =|a| Y_s, \quad \nabla_s Y_s=-|a| X_s
$
where $X_s$ and $Y_s$ are unit vector fields, $X_s$ is the tangent vector
to the curve $\Gamma$ parametrised by the arc-length $s$, $\nabla_s$ 
is the covariant derivative of $g$ along $\Gamma$ and the positive
constant $|a|^{-1}$ is the {\em radius} of the circle $\Gamma$.

If $\Gamma^s_0$ is the parallel displacement of tangent vectors along
$\Gamma$ from $\Gamma(s)$ to $\Gamma(0)$, and
${X_s}^*=\Gamma^s_0(X_s)$, then the development $\Gamma^*$ is the unique
curve in $T_x M$ starting at the origin such that its tangent vector is
parallel to $X_s$ in the Euclidean sense.
} 
of a circle $\Gamma$ starting at a point $x\in M$ is
an ordinary circle in $T_xM$.
It is also the case that circles on the round 
sphere $S^n$ are intersections of $S^n\subset \R^{n+1}$ with planes 
(not necessarily through the origin) in $\R^{n+1}$. Thus they are indeed
circles.
\subsection{First Integrals}
If $Y\in \Lambda^2(M)$ satisfies the conformal Killing--Yano (CKY) equation
\be
\label{CKY}
\nabla_aY_{bc}=\nabla_{[a}Y_{bc]}-2g_{a[b}K_{c]}
\ee
for some one--form $K\in \Lambda^1(M)$, then (if $\mbox{dim}(M)=4$) equation (\ref{conf_circ}) implies that
\be
\label{Q}
Q=u^aa^bY_{ab}-u^aK_a
\ee
is constant along the conformal geodesics. In the special case where 
$Y_{ab}$ is a K\"ahler form, the one form $K$ is zero and the linear
term in $(\ref{Q})$ is not present. In this case the first integral
$(\ref{Q})$ has been called the {\em complex torsion} (though it is real)  in \cite{jap2}.
In four dimensions the condition $u^c\nabla_c Q=0$ with $Q$ given by 
(\ref{Q}) was established in \cite{tod_circles}.
In \cite{rodar} it was put in the general context of parabolic geometries.
\subsection{Examples}
\subsubsection{Circles on $\R^n$.} These are just ordinary circles. 
Equation (\ref{conf_circ_e}) becomes $\dddot{{\bf x}}=-|a|^2\dot{\bf x}$,
and can be readily solved
\[
{\bf{x}}(s)={\bf x}(0)+|a|^{-1}{\bf v}(0)\sin{(|a|s)}
+|a|^{-2}{\bf a}(0)(1-\cos{(|a|s)}),
\]
where ${\bf x}(0)$ is arbitrary,  ${\bf v}(0)$ is a unit vector in $\R^n$, and the vector 
${\bf a}(0)$ is orthogonal to ${\bf v}(0)$, and has squared norm $|a|^2$.
We recognise these conformal geodesics as circles centered at
${\bf x}(0)+|a|^{-2}{\bf a}(0)$, and with Euclidean radius $|a|^{-1}$.

Conformal rescalings of the metric preserve conformal geodesics, so circles
on $S^n$ are also circles (intersections of $S^n\subset \R^{n+1}$ with planes 
in $\R^{n+1}$ which do not necessarily pass through the origin). In particular
all conformal geodesics on $S^n$ are closed.
\subsubsection{Magnetic motion on the upper half--plane.}
\label{cccm} 
Conformal geodesics on the hyperbolic space can also be obtained from those on 
$\R^n$. This time however they are not necessarily closed, and we shall
consider this case separately to illustrate the significance of the first 
integrals in the qualitative behaviour of conformal geodesics.

Let ${\mathbb{H}}$ be the upper half--plane with  coordinates  $(x, y)$, and a constant curvature
metric
\[
g=\frac{dx^2+dy^2}{y^2}, \quad y>0
\]
(so the Ricci scalar equal to $-2$) and a magnetic field
$F=B\; \mbox{vol}_{\mathbb{H}}$ given by a constant $B>0$ multiple of
the parallel volume form on ${\mathbb{H}}$. We chose the potential
$\Phi=By^{-1}dx$, so that $F=d\Phi$. The magnetic 
Lagrangian\footnote{The variational formulation of the conformal geodesic equations with $n>2$ has been developed
in \cite{DK21}.}
\[
L=\frac{1}{2}\frac{\dot{x}^2+\dot{y}^2}{y^2}-B\frac{\dot{x}}{y}
\]
admits a conserved energy integral which we set to $1$.
Taking $u$ to be a unit tangent vector to an integral curve
of the corresponding Euler--Lagrange equations we find
\[
a\equiv\nabla_u u= B J(u),
\]
where $J$ is the complex structure on ${\mathbb{H}}$ defined by
$g(v, J(w))=\mbox{vol}_{\mathbb{H}}(v, w)$. Therefore the circle equations
(\ref{conf_circ_e}) hold with $|a|^2=B^2$, as
\[
\nabla_u a=B \nabla_u J(u)=B J(\nabla_u u)=B^2 J^2(u)=-B^2 u.
\]
The two first
integrals $|a|^2$, and $Q$  given by (\ref{Q})  (where now $Y=F$) are therefore equal as
\[
Q\equiv F(u, a)=B\; \mbox{vol}_{\mathbb{H}}(u, B J(u))=B^2 g(u, u)=B^2.
\]
The resulting trajectories are circles of radius $|B|^{-1}$.
The behaviour of the circles depends
on the value of $B$ (Figure 1)
\begin{center}
\includegraphics[width=7cm,height=6cm,angle=0]{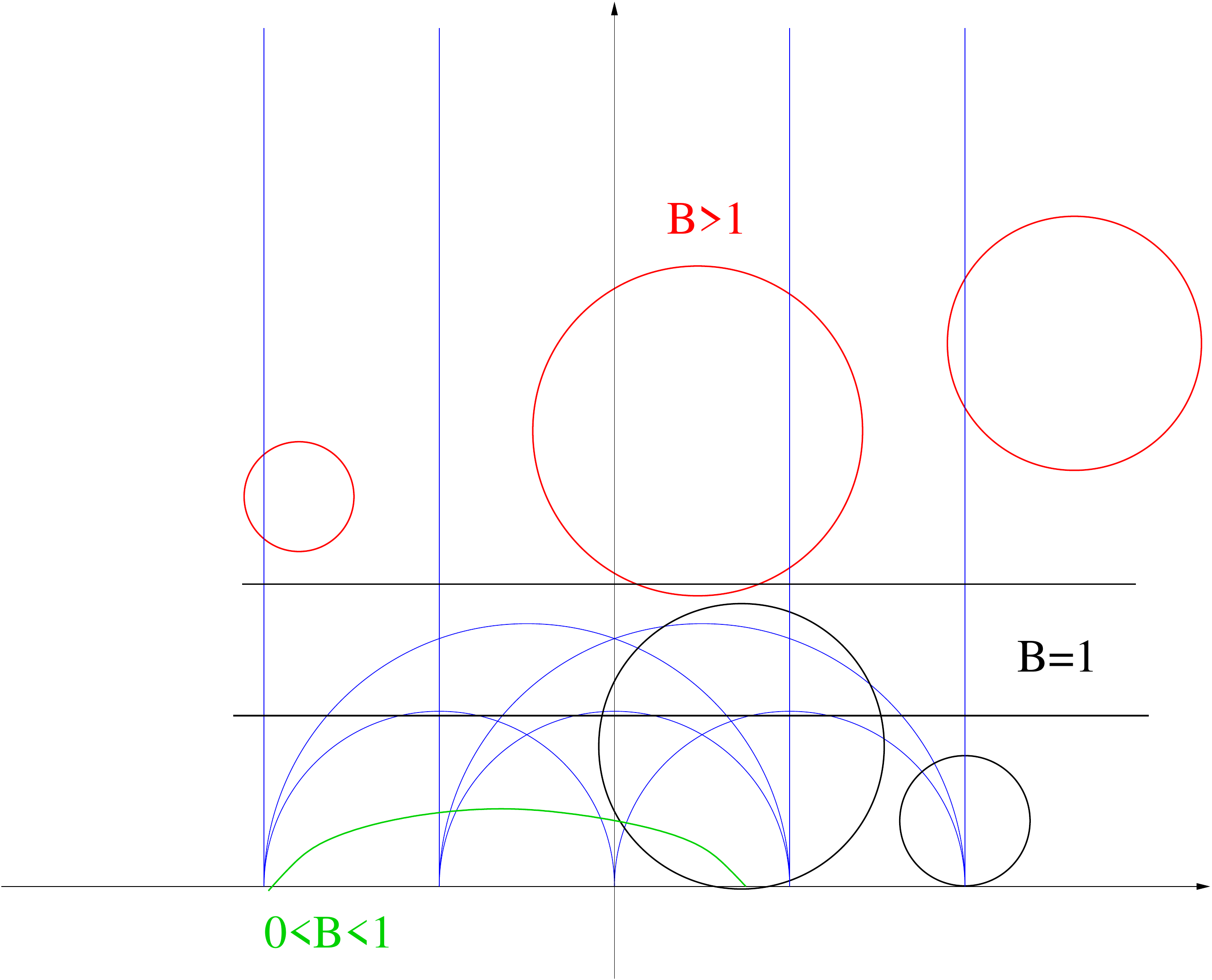}
\begin{center}
{\em Figure 1.} Geodesics (in blue) and circles on ${\mathbb{H}}$.
\end{center}
\end{center}
The circles are open (in the language
of \cite{comet} the magnetic field is not strong enough to capture the 
particle) and unbounded if $0<B<1$, and closed if $B>1$. The special case
$B=1$ corresponds to horocircles tangent to the $y=0$ boundary
of  ${\mathbb{H}}$, or lines $y=$const parallel to the $x$--axis.
\section{Circles on hyper--K\"ahler four--manifolds and the Lorentz 
force}
Let $(M, g)$ be an oriented Riemannian--four manifold with a hyper--K\"ahler
structure given by parallel two--forms
$\Omega^1, \Omega^2, \Omega^3$. In particular $(M, g)$ is Ricci--flat,
so the conformal geodesic equations reduce to (\ref{conf_circ_e}).
\begin{prop}
\label{prop31}
The conformal geodesic equations on a hyper--K\"ahler four--manifold
$(M, g)$ reduce to the Lorentz force equations
\be 
\label{lorentzf}
\nabla_u u=|a|J(u),
\ee
where $J$ is the complex structure corresponding to the K\"ahler form given by
a constant self--dual Maxwell field
$F=\sum_j c^j\Omega^j$ with  constants
$(c^1, c^2, c^3)$ such that $|{\bf c}|=1$. 
\end{prop}
\noindent
{\bf Proof.}
Each of the K\"ahler forms gives rise to a first
integral (\ref{Q}), so there are three of those.
These three integrals allow us to integrate the circle equations once, and
reduce them, on a level set 
\be
\label{thecs}
c^i=\Omega^i(u, a), \quad i=1, 2, 3,
\ee
to geodesic equations of a charged 
particle moving in a magnetic field.

If the orientation of $(M, g)$ is chosen so that the Riemann tensor of $g$ is
anti--self--dual, then this magnetic field is self--dual. To see how the Lorentz
force arises assume that $M$ is paralelizable  chose an orthonormal frame such that
\[
g={(e^1)}^2+{(e^2)}^2+{(e^3)}^2+{(e^4)}^2,
\]
and
\[
\Omega^1=e^1\wedge e^4+e^2\wedge e^3, \quad
\Omega^2=e^2\wedge e^4+e^3\wedge e^1, \quad
\Omega^3=e^3\wedge e^4+e^1\wedge e^2.
\]
Let $u^a$ and $a^a$ be the components of $u, a$ in the frame of vector--fields
dual to $e^a$. The relations (\ref{thecs}) can then be solved for $a$.
Using the ordinary vector notation for the components
${\bf u}=(u^1, u^2, u^3), {\bf a}=(a^1, a^2, a^3)$ we can rewrite  (\ref{thecs}) 
as  ${\bf c}={\bf u}a^4-{\bf a}u^4+{\bf u}\wedge {\bf a}$, and find
\be
\label{aGH}
{\bf a}= {\bf c}\wedge {\bf u} - u^4{\bf c} , \quad a^4={\bf u}\cdot {\bf c}.
\ee
Note that
\be
\label{canda}
g(a, u)={\bf a}\cdot {\bf u}+a^4u^4=0, \quad g(a, a)=|{\bf c}|^2
\ee
so the first integrals $g(a, a)$ and ${\bf c}$ are not independent. The circle equations now reduce to $a=\nabla_u u$, or (\ref{lorentzf})
which is the Lorentz force equation for a unit mass particle with charge
$e=|a|$ in a constant self-dual Maxwell field
\[
F=\frac{1}{|a|}\Big(c^1\Omega^1+c^2\Omega^2+c^3\Omega^3).
\]
In view of (\ref{canda}) we can redefine ${\bf c}$ to be a unit--vector,
and recover the statement of Proposition \ref{prop31}.
\koniec
For an uncharged particle $|a|=0$, and (\ref{lorentzf})  reduces to the ordinary geodesic equation.
To make further progress with $|a|\neq 0$ we need to
seek first integrals of (\ref{lorentzf}). Any Killing vector
on $(M, g)$ which is also tri--holomorphic, i .e. 
${\mathcal L}_K\Omega^i=0$, will give rise to a first integral: If
${F}=d\Phi$, the potential one--form $\Phi$ can be chosen such that
${\mathcal L}_K\Phi=0$ (we then say that $K$ Lie drags $\Phi$). Indeed, computing
\[
0={\mathcal L}_K F= d(K\hook F) \quad \mbox{implies that, locally}\quad K\hook F=df
\]
for some function $f:M\rightarrow \R$. Now
\[
{\mathcal L}_K\Phi=d(f+K\hook \Phi).
\] 
Choosing $f_1:M\rightarrow \R$ such that $K(f_1)=-f-K\hook \Phi$,  and performing a gauge transformation 
$\Phi\rightarrow \Phi+df_1$ we establish that in this gauge ${\mathcal L}_K\Phi=0$. This implies that
\be
\label{cg4}{\mathcal K}:=K^a(u_a+e\Phi_a)=K^ap_a
\ee
is conserved along conformal geodesics on a level set of the $c^i$. The momentum $p_a$ is defined by 
\[u^a=g^{ab}(p_a-e\Phi_a),\]
which is the usual convention for charged particles. 

For (\ref{lorentzf}) to be integrable we in particular need the geodesic equation to be integrable, and we need four constants of the motion in involution. One is the 
Hamiltonian (the
unit mass) and the other three need to come from Killing vectors and Killing tensors. If at least two of them arise from Killing vectors then, to be involutive, these must commute and to give integrals for (\ref{lorentzf}) they must Lie-drag $\Omega^i$. Now necessarily there is a constant linear combination of them which has purely anti-self-dual covariant derivative (this is Theorem 2.3 in 
\cite{gibbons_ruback}, credited to Hitchin by these authors) and this is the characterising property of the Gibbons-Hawking 
metrics \cite{GHr}. Consequently we henceforth restrict to these metrics, but there may be other integrable cases with two or more Killing tensors (see  \cite{Houri, frolov} for other applications for  Yano tensors to geodesic integrability). 
\subsection{The Gibbons Hawking metrics}
All
hyper--K\"ahler four--manifolds which admit  a tri--holomorphic Killing vector can be put
in the Gibbons--Hawking form \cite{GHr}
\be
\label{GH}
g=V {\bf dx}\cdot {\bf dx}+V^{-1}(d\tau+ \omega)^2
\ee
where $(V, \omega)$ is a function, and a one--form on $\R^3$. Choosing a 
frame
\[
e^4=V^{-1/2}(d\tau+\omega), \;\; e ^i=V^{1/2}dx^i, \quad E_4=V^{1/2}\p_\tau, \;\; E_i= V^{-1/2} (\p_i-A_i\p_\tau) 
\]
one verifies that the two--forms ${\Omega^i}$ are closed iff
\[
d\omega=*dV,
\]
which in particular implies that $V$ is harmonic.

Consider the Killing vector $K=\partial_\tau$ of metric (\ref{GH}), and set $K^\flat=V^{-1}(d\tau+\omega)$ and then its derivative is
anti-self-dual (the characterising property of the Gibbons-Hawking metrics). It follows that $K$ is tri--holomorphic, i. e. 
\[{\mathcal{L}}_K\Omega^i=0,\] where $\mathcal{L}_X$ is the Lie-derivative along the vector field $X$.
It therefore follows that $K$ gives rise to three moment maps: for each $i$ there is a function $x^i$ satisfying
\be\label{h1}
K\hook \Omega^i=d x^i
\ee
and these $x^i$ are (up to scale) the flat coordinates $(x,y,z)$ in (\ref{GH}).
\vskip5pt
The Killing vector $K$ will give a constant 
of the motion (\ref{cg4}) for (\ref{lorentzf}).
This is conserved for any $V$, and we will return to it when we consider examples below. It is however not sufficient for 
integrability. 

To integrate (\ref{lorentzf}) we need either two more Killing vectors, commuting with each other and with $K=\partial_\tau$, or another Killing vector commuting with $K$ and a Killing tensor with suitably involutive. The first case is not interesting (see the Appendix A) but for the second we can take motivation from the familiar fact that both the self-dual Taub-NUT metric and the Eguchi-Hanson metric admit Killing tensors \cite{gibbons_ruback} to concentrate on these two metrics in
 \S \ref{sectionTN}  and \S \ref{sectionEH}.

\section{Circles on anti--self--dual Taub--NUT}
\label{sectionTN}
The ASD Taub--NUT metric \cite{GHr} 
is a special case of (\ref{GH}) with $V=1+m/r$ where $m$ is a non--negative constant, and (using spherical
polar coordinates on $\R^3$) $\omega=m\cos{\theta}  d\phi$ and

\[
g=\Big(1+\frac{m}{r}\Big)(dr^2+r^2(d\theta^2+\sin^2{\theta} d\phi^2))+ \Big(1+\frac{m}{r}\Big)^{-1}(d\psi+m\cos{\theta} d\phi)^2
\]
writing $\psi$ for $\tau$ to emphasise the Bianchi IX $SU(2)$ symmetry.
The singularity at $r=0$ is removable. The infinity $r\rightarrow \infty$ has the topology of a one--monopole $S^1$ bundle over $S^2$. This behaviour is referred to as
asymptotic local flatness (ALF). Thus ASD Taub--NUT is regular everywhere, and is an example of a gravitational instanton.

The metric has $SU(2)\times U(1)$ as its group of 
isometries
where the $U(1)$ action generated by $\p/\p \psi$ is tri--holomorphic, i. e.
\[
{\mathcal{L}}_{\p/\p\psi} \Omega^i=0,\quad i=1, 2, 3
\]
and so the self--dual derivative of the corresponding one--form 
$g(\p/\p\psi, \cdot)=V^{-1}(d\psi+\cos{\theta}d\phi)$ vanishes. The $SU(2)$ action
is not tri--holomorphic, and the exterior derivatives of the corresponding Killing vectors
do not have a definite duality (but their SD derivatives are parallel). However for any unit vector
${\bf c}\in \R^3$ there exists a linear combination $L$ of the isometric generators of $SU(2)$, such
that $L$ Lie--drags $\sum_ic^i\Omega^i$.
\subsection{Conformal Killing--Yano tensors and Arnold--Liouville integrability} 
The ASD Taub--NUT metric admits a Killing--Yano $2$--form $Z$, i. e. a
two--form such that
\[
\nabla_a Z_{bc}=\nabla_{[a}Z_{bc]}.
\]
The self--dual, and anti--self--dual parts of $Z$ are conformal Killing--Yano
two--forms which satisfy (\ref{CKY}):
\begin{eqnarray*}
Z&=& (d\psi+m\cos{\theta} d\phi) \wedge dr+(2r+m)(r+m)r\sigma_1\wedge\sigma_2\\
&=& Y-W, \quad\mbox{where}\qquad *Y=Y, \quad *W=-W 
\end{eqnarray*}
and
\be
\label{Wtensor}
Y=\frac{r^3}{m} d(V(d\psi+m\cos{\theta} d\phi)),\quad W=-\frac{(r+m)^3}{m} d(V^{-1}(d\psi+m\cos{\theta} d\phi)).
\ee
The left invariant one--forms $\sigma_i, i=1, 2, 3$ 
on the group manifold $SU(2)$ are such that
\be
\label{cijk}
d\sigma_1+\sigma_2\wedge\sigma_3=0, \quad
d\sigma_2+\sigma_3\wedge\sigma_1=0, \quad d\sigma_3+\sigma_1\wedge\sigma_2=0.
\ee
These one--forms can be represented in terms of Euler angles
by
\[
\sigma_1+i\sigma_2=e^{-i\psi}(d\theta+i\sin{\theta} d\phi), \qquad
\sigma_3=d\psi+\cos{\theta}d\phi,
\]
where to cover $SU(2)=S^3$ we require the ranges
$
0\leq\theta\leq\pi, \quad 0\leq\phi\leq 2\pi, \quad
0\leq\psi\leq 4\pi.$
Up to a multiple of $4$ the two--form $Z$ coincides with the Killing--Yano 
two--from found by Gibbons and Ruback \cite{gibbons_ruback}
(see formula (2.12) in this reference with $n=1$, and make
a coordinate transformation $t=m+2r$. See also \cite{valent}).

 Both $Y$ and $W$ satisfy (\ref{CKY}) with the same Killing vector $K=g(\p/\p\psi, \cdot)$.
The CKY two--form $Y$ was discovered in \cite{DTeinstein}, where it was used to show that the
ASD Taub--NUT metric is conformal to a scalar--flat K\"ahler metric with a 
non--constant conformal factor. Computing the first integral (\ref{Q}) corresponding to $Y$ gives an expression which is functionally dependent with the first
integrals arising from the Killing vectors $K=\p/\p\psi$ and $L=\p/\p\phi$, so it does 
not give anything useful. We instead turn to $W$, where the resulting first 
integal (\ref{Q}) is non--trivial, but it does not in general commute with the 
first integrals (\ref{cg4}) arising from $\p_\psi$ and $\p_\phi$. To get around this problem we proceed as follows:

Given (\ref{lorentzf}) with some $F$, we may first rotate the $(x, y, z)$ coordinates until $F=-\Omega^3$: if this case is integrable then the general case is integrable. Introducing spherical polar coordinates
$(r, \theta, \phi)$ in place of $(x, y, z)$ yields
\begin{eqnarray}
F&=& -\Omega^3=d\Phi\quad\mbox{with}\nonumber\\
\Phi&=&-z(d\psi+\omega)+\Big(\frac{m}{r}+\frac{1}{2}\Big)(ydx-xdy)\nonumber\\
&=&-r\cos{\theta}d\psi-\Big(mr+\frac{1}{2}r^2\sin^2{\theta}\Big)d\phi
\end{eqnarray}
or in polars and with the index raised
\[\Phi^\#=-zV\partial_\psi+\frac{1}{V}\partial_\phi.\]
In this form it is clear that $\Phi^\#$ commutes with the Killing vectors $K:=\partial_\psi$ and $L:=\partial_\phi$, which in turn commute with each other, leading to two constants of the motion for (\ref{lorentzf}). 

This  leads to four independent first integrals 
${\mathcal I}_a=({\mathcal K}, {\mathcal L}, {\mathcal W}, {\mathcal H})$
two of which are linear in the velocity, and two of which are quadratic:
\begin{eqnarray}
\label{first_ints}
{\mathcal K}&=&K^a(u_a+e\Phi_a), \quad {\mathcal L}=L^a(u_a+e\Phi_a), \\
{\mathcal H}&=&\frac{1}{2}g_{ab}u^a u^b, \quad
{\mathcal W}= eW_{ac}{{F}^c}_b u^a u^b-2{\mathcal H}\;K^a u_a.\nonumber
\end{eqnarray}
To complete the calculation we shall use  Hamiltonian formalism. In the presence of a
magnetic field there are two ways to proceed starting from the magnetic Lagrangian
${\mathcal G}=(1/2) g_{ab} u^a u^b+e\Phi_b u^b$. One way is to define the conjugate momentum 
\[
p_a=\frac{\p {\mathcal G}}{\p u^a}=u_a+e\Phi_a
\]
and perform a Legendre transform to find four first  integrals
\begin{eqnarray*}
{\mathcal K}&=&p_\psi,  \quad {\mathcal L}=p_\phi, \\
 {\mathcal H}&=&\frac{1}{2}g^{ab}(p_a-e\Phi_a)(p_b-e\Phi_b), \quad
{\mathcal W}=e{W^a}_c F^{cb} (p_a-e\Phi_a)(p_b-e\Phi_b)-2{\mathcal H}({\mathcal K} -e K^a\Phi_a).
\end{eqnarray*}
These are in involution with respect to the standard symplectic structure on $T^*M$ given by $dp_a\wedge dq^a$, where
$q^a=(r, \phi, \theta, \psi)$ and $p_a=(p_r, p_\phi, p_\theta, p_\psi)$.
Alternatively define a momentum by
\[
P_a=p_a-e\Phi_a,
\]
and consider the charged symplectic form
\be
\label{charged}
{\bf \Omega}= dP_a\wedge dq^a+eF
\ee
on the $8$--dimensional phase--space with coordinates
$q^a=(r, \phi, \theta, \psi)$ and $P_a=(P_r, P_\phi, P_\theta, P_\psi)$.
The set of first integrals is then ${\mathcal I}_a=({\mathcal K},  {\mathcal L}, {\mathcal H}, {\mathcal W})$, where
\begin{subequations}
\label{mag_integrals}
\begin{eqnarray}
{\mathcal K}&=&P_\psi-e r\cos{\theta}, \\
 {\mathcal L}&=&P_\phi-e \Big(mr+\frac{1}{2}r^2\sin^2{\theta}\Big),\\
 {\mathcal H}&=&\frac{1}{2}\Big(\frac{r}{r+m} {P_r}^2\nonumber\\
& &+ \frac{1}{r(r+m)}{P_\theta}^2+\frac{r+m}{r}{P_\psi}^2
+\frac{1}{r(r+m)\sin^2{\theta}}(P_\phi-m\cos{\theta} P_\psi)^2\Big),\\
{\mathcal W}&=&\frac{e\cos{\theta}}{r}(rP_r-\tan{\theta} P_\theta)^2-\frac{e}{r\cos{\theta}}{P_\theta}^2\nonumber\\
& &
-\frac{e\cos{\theta}}{r\sin^2{\theta}}({P_\phi}^2+(m^2-r^2\sin^2{\theta}){P_{\psi}}^2-2\frac{m+r\sin^2{\theta}}{\cos{\theta}}P_\psi P_\phi\Big)-2{\mathcal H}P_\psi. 
\end{eqnarray}
\end{subequations}
These integrals are in involution
\[
\{ {\mathcal I}_a,  {\mathcal I}_b
\}=0, \quad a, b =1 ,\dots, 4
\]
with respect to Poisson brackets of the charged symplectic structure 
(\ref{charged}). We have established
\begin{prop}
\label{prop_TN1}
The conformal geodesic equation on the anti--self--dual Taub--NUT manifold is completely integrable
in the sense of Arnold--Liouville.
\end{prop}
 Rather than finding the action--angle variables  explicitly, we follow \cite{gibbons_ruback} and seek to separate the Hamilton-Jacobi equation, but for (\ref{lorentzf}) rather than for the geodesic equation.
\subsection{Separating the Hamilton-Jacobi equation for anti-self-dual Taub-NUT}
Following Gibbons-Ruback \cite{gibbons_ruback}, we introduce parabolic coordinates
\[\eta=r(1+\cos\theta),\;\;\xi=r(1-\cos\theta),\;\;\phi=\phi,\]
so that in the cylindrical polars obtained from $(x,y,z)$ in (\ref{GH})
\[z=\frac12(\eta-\xi),\;\;\rho=\sqrt{\eta\xi},\;\;\phi=\phi.\]
Now
\[z+i\rho=\frac12(\sqrt{\eta}+i\sqrt{\xi})^2,\]
and
\[dz^2+d\rho^2=\frac14(\eta+\xi)\big(\frac{d\eta^2}{\eta}+\frac{d\xi^2}{\xi}\big),\]
\[V=1+\frac{m}{r}=\frac{(\eta+\xi+2m)}{(\eta+\xi)},\quad\omega=m\cos\theta d\phi=m\frac{(\eta-\xi)}{(\eta+\xi)}d\phi,\]
so that the anti-self-dual Taub-NUT metric  is
\[g=V\left(\frac14(\eta+\xi)\big(\frac{d\eta^2}{\eta}+\frac{d\xi^2}{\xi}\big)+\eta\xi d\phi^2\right)+\frac{1}{V}\left(d\psi+m\frac{(\eta-\xi)}{(\eta+\xi)}d\phi\right)^2.\]
For the geodesic equation first, we seek Hamilton's Principal Function $S$ in the separated form
\[S=E\psi+J\phi+F(\xi)+G(\eta),\]
with constant $E,J$ then the Hamilton-Jacobi equation is
\[\mu^2=g^{ab}\nabla_aS\nabla_bS
\]
where $\mu$ is the conserved mass which we will always take to be one.

Algebra leads to
\[4\big(\eta G_\eta^2+\xi F_\xi^2)=\mu^2(\eta+\xi+2m)-J^2(\frac{1}{\eta}+\frac{1}{\xi})+2JEm(-\frac{1}{\eta}+\frac{1}{\xi})-E^2(\eta+\frac{m^2}{\eta}+\xi+\frac{m^2}{\xi}+4m),\]
and, confirming \cite{gibbons_ruback}, this separates with a new constant $Q$:
\[4\eta G_\eta^2=Q+\mu^2(\eta+m)-\frac{J^2}{\eta}-\frac{2JEm}{\eta}-E^2(\eta+\frac{m^2}{\eta}+2m),\]
\[4\xi F_\xi^2=-Q+\mu^2(\xi+m)-\frac{J^2}{\xi}+\frac{2JEm}{\xi}-E^2(\xi+\frac{m^2}{\xi}+2m).\]
The constant $Q$ is quadratic in momenta and must be associated with a quadratic Killing tensor. To obtain the geodesics from $S$, we solve
\[\dot{q}^a=g^{ab}p_b=g^{ab}\nabla_b S,\]
so that
\[\dot\eta=\frac{4\eta}{\eta+\xi+2m}G_\eta,\;\;\dot\xi=\frac{4\xi}{\eta+\xi+2m}F_\xi,\;\;p_\phi=J,\;\;p_\psi=E.\]
We shall now modify this separation procedure to include the Maxwell field
present for conformal geodesic motion.
\begin{prop}
\label{prop_TN2}
The Hamilton--Jacobi equation for the conformal geodesics on anti--self--dual Taub-NUT manifold is separable.
\end{prop}
\noindent
{\bf Proof.}
To include a Maxwell field, 
we follow the Hamilton--Jacobi formalism in the presence of electro-magnetic field (see e.g. \cite{LL}). The method is to suppose the that potential is $\Phi_a$, so that
\[ \dot{q}^a=g^{ab}(p_b-e\Phi_b)\quad\mbox{and}\quad 
H=\frac12g^{ab}(p_a-e\Phi_a)(p_b-e\Phi_b).\]
Now set $p_a=\nabla_a S$ and take the Hamilton-Jacobi equation to be
\[\mu^2=g^{ab}(\nabla_a S-e\Phi_a)(\nabla_b S-e\Phi_b).\]
For us
\[\Phi=-z(d\psi+\omega)-r^2\sin^2\theta\big(\frac12+\frac{m}{r}\big)d\phi=-\frac12(\eta-\xi)d\psi-\frac12(\eta\xi+m(\eta+\xi))d\phi,\]
\[=Xd\phi+Yd\psi, \mbox{  say.}\]
The changes to the Hamilton-Jacobi equation are
\[J\rightarrow J-eX,\;\;E\rightarrow E-eY,\]
leading to
\begin{eqnarray}
\label{sepwithe}
&& 4\big(\eta G_\eta^2+\xi F_\xi^2)=\mu^2(\eta+\xi+2m)-J^2(\frac{1}{\eta}+\frac{1}{\xi})+2JEm(-\frac{1}{\eta}+\frac{1}{\xi})-E^2(\eta+\frac{m^2}{\eta}+\xi+\frac{m^2}{\xi}+4m)\nonumber\\
&&-eJ(\eta+\xi+4m)-eE(\eta^2-\xi^2+3m(\eta-\xi))-\frac{e^2}{4}(\eta(\eta+2m)^2+\xi(\xi+2m)^2),
\end{eqnarray}
which still separates with a quadratic constant and implied Killing tensor. As before, for the particle paths solve
\[
\dot{q}^a=g^{ab}p_b=g^{ab}\nabla_b S.
\]
\koniec
Thus in anti-self-dual Taub-NUT the conformal geodesic equation is completely integrable: it admits four first integrals
in involution established in Proposition \ref{prop_TN1}, and the associated Hamilton--Jacobi equation is separable which reduces
the integration of the conformal geodesic equations to quadratures. This completes the proof of Theorem \ref{thm_TN}.
\vskip 5pt
Separating (\ref{sepwithe}) with a constant $Q$ gives
\be
\label{sep123}
\dot{\eta}=\frac{4\eta G_\eta}{\eta+\zeta+2m}, \quad \dot{\zeta}=\frac{4\zeta F_\zeta}{\eta+\zeta+2m}
\ee
with
\[
4\eta G_\eta= U(\eta, E, Q), \quad 4 \zeta F_\zeta= U(\zeta, -E, -Q)
\]
where
\begin{eqnarray}
\label{Usep}
U(x, E, Q)&=&\sqrt{u_0+u_1 x+ u_2 x^2+u_3 x^3+u_4 x^4}, \\
u_0&=&-4(J+Em)^2, \quad u_1=4(Q+m(\mu^2+2Je-2E^2)),\nonumber\\ 
u_2&=&4(\mu^2-m^2e^2-E^2+Je+3Eme), \quad
u_3=4e(E-me), \quad u_4=-e^2.\nonumber
\end{eqnarray}
Therefore (\ref{sep123}) imples that the unparametrised conformal geodesic equations are solvable in terms
of elliptic functions:
\[
\int\frac{d\eta}{U(\eta, E, Q)}=\int\frac{d\zeta}{U(\zeta, -E, -Q)}.
\]
\section{Circles on Eguchi--Hanson}
\label{sectionEH}
The Eguchi--Hanson (EH) metric \cite{EHref}
\[
g=f^{-2}dr^2+\frac{1}{4}r^2(\sigma_1^2+\sigma_2^2)+\frac{1}{4}r^2f^2\sigma_3^2,\quad \mbox{where}\quad f=\sqrt{1-\frac{\alpha^4}{r^4}}, \quad \alpha=\mbox{const}
\]
is another special case of (\ref{GH}).  It corresponds to 
\be
\label{VEH}
V=\frac{1}{r_1}+\frac{1}{r_2},
\ee
where $r_1, r_2$ are Euclidean distances from points $P_1, P_2$
which we can chose to be
\[
P_1=(0, 0, -\alpha), \quad P_2=(0, 0, \alpha)
\]
for constant positive $\alpha$. To achieve regularity the period of the $\psi$ coordinate
in $\sigma_3$ should be $2\pi$ rather than $4\pi$. Therefore the surfaces
of constant $r$ are real projective spaces, and at large $r$ the metric looks like
$\R^4/\Z_2$ rather than  Euclidean space. The Eguchi--Hanson is an example of an
asymptotically locally Euclidean (ALE) manifold (see e.g. Chapter 9 of \cite{Dbook}).

The isometry group of EH is  enhanced to
$SO(3)\times U(1)$, but this time it is the $SO(3)$ which acts tri--holomorphically (so that
EH can be put in the form (\ref{GH}) in many ways). The $U(1)$ action is not tri--holomorphic.
The parallel basis of ${\Lambda^2}_+$ is given by 
\[
\Omega^{i}=e^i\wedge e^4+\frac{1}{2} \varepsilon^{ijk} e^j\wedge e^k,
\]
with
\be
\label{basisEH}
e^1=\frac{1}{2}r\sigma_1, \quad e^2=\frac{1}{2}r \sigma_2, \quad
e^3=\frac{1}{2} rf\sigma_3,\quad  e^4=f^{-1}dr.
\ee
\begin{prop}
\label{propEH}
All conformal geodesics on the $SO(3)$ orbits of constant $r$ in the Eguchi--Hanson manifold are of the form
\begin{eqnarray}
\label{starstar}
&&{(u^1)}^2+{(u^2)}^2+{(u^3)}^2=1, \quad
c_1 u^1+c_2 u^2=\frac{rf}{2(1-f^2)} (h^2-2({c_1}^2+{c_2}^2)-2m_1),\nonumber\\
&&\mbox{where}\quad
h=\frac{1-f^2}{fr}u^3
+\frac{m_0}{{c_1}^2+{c_2}^2}  
\end{eqnarray}
and $(u^1, u^2, u^3)$ depend on $(\phi, \psi, \theta)$ as in (\ref{uspaul})
where
$(c_1, c_2, c_3, m_1, m_2, p_1, p_2, p_3)$ are constants of integration.
\end{prop}
We shall first perform the general analysis of the circle equations
using the Lorentz force system (\ref{lorentzf}), and reduce the equations
to a system of three 1st order ODEs (\ref{finaleq}) for three unknown functions. We 
shall then show that this system is completely solvable under the additional 
assumption $r$=const. This will give the proof of Proposition \ref{propEH}.

With $a$ determined in terms of $u$ by (\ref{aGH})  (note that now the components of $a$ and $u$ refer to the basis \ref{basisEH})
the circle equations reduce to the definition of the acceleration
$a=\nabla_u u$, which 
becomes
\begin{subequations}
\label{ueq1}
\begin{eqnarray}
\dot{u}^1&=&  2\frac{k^4}{ r\sqrt{1-k^4}}u^2u^3
-\frac{\sqrt{1-k^4}}{r}u^1u^4 +c^2u^3-c^3u^2-c^1u^4,\\
\dot{u}^2&=&   -2\frac{k^4}{ r\sqrt{1-k^4}}u^1u^3-
\frac{\sqrt{1-k^4}}{r}u^2u^4 +c^3u^1-c^1u^3-c^2u^4,\\
\dot{u}^3&=& -\frac{1+k^4}{r\sqrt{1-k^4}}u^3u^4 +c^1u^2-c^2u^1-c^3u^4,\\
\dot{u}^4&=&\frac{\sqrt{1-k^4}}{r}(1-(u^4)^2)+
2\frac{k^4}{ r\sqrt{1-k^4}}(u^3)^2+c^1u^1+c^2u^2+c^3u^3
\end{eqnarray}
\end{subequations}
where $k\equiv\alpha/r$.
In this case we can also find the potentials of the triple of parallel two--forms
\begin{subequations}
\label{potentials}
\begin{eqnarray}
\Omega^1&=&  e^1\wedge e^4+e^2\wedge e^3=d\Big(-\frac{rf}{2}e^1\Big)\\
\Omega^2&=&  e^2\wedge e^4+e^3\wedge e^1=d\Big(-\frac{rf}{2}e^2\Big)\\
\Omega^3&=&  e^3\wedge e^4+e^1\wedge e^2=d\Big(-\frac{r}{2f}e^3\Big).
\end{eqnarray}
\end{subequations}
The relation
\[
u=\dot{\phi}\frac{\p}{\p\phi}+\dot{\theta}\frac{\p}{\p\theta}+
\dot{\psi}\frac{\p}{\p\psi}+\dot{r}\frac{\p}{\p r}=u^aE_a
\]
gives
\begin{subequations}
\label{uscor}
\begin{eqnarray}
u^1&=&\frac{1}{2}r(\sin{(\theta)}\sin{(\psi)}\dot{\phi}+\cos{(\psi)}\dot{\theta})\\
u^2&=&\frac{1}{2}r(\sin{(\theta)}\cos{(\psi)}\dot{(\phi)}-\sin{(\psi)}\dot{\theta})\\
u^3&=&\frac{1}{2}rf(\dot{\psi}+\cos{(\theta)}\dot{\phi})\\
u^4&=&\frac{1}{f}\dot{r}.
\end{eqnarray}
\end{subequations}
\subsection{Magnetic first integrals, and right invariant vector fields}
The equations (\ref{ueq1}) correspond to a forced geodesic motion 
of a charged particle in a constant
self--dual  Maxwell field $F=c^1\Omega^1+c^2\Omega^2+c^3\Omega^3$
\be
\label{lorenzmd}
u^b\nabla_b u^a= {F^a}_b u^b.
\ee
Let $\Phi$ be a potential for $F=d\Phi$. Using (\ref{potentials}) we
find
\[
\Phi=-\frac{1}{2}\Big(c_1rfe^1+c_2rfe^2+c_3\frac{r}{f}e^3\Big).
\]
For any tri--holomrphic Killing vector $R$, such that 
${\mathcal{L}}_R \Phi=0$ (this can always be arranged by adding a total derivative of a function to $\Phi$), the function
\be
\label{pauls_c}
p\equiv R^a(u_a+\kappa  \Phi_a)
\ee
is a first integral of (\ref{lorenzmd}) (here $\kappa$ is a numerical factor
and the value
of $\kappa$ can be absorbed into constants $c^i$).

Given the left--invariant one--forms $\sigma_i$, let $L_i$ be the dual basis
of left--invariant vector fields, i. e. $L_i\hook \sigma_j=\delta_{ij}$, and
\[
[L_i, L_j]=\frac{1}{2}\epsilon_{ijk}L_k.
\]
The tri--holomorphic and isometric $SO(3)$ action is generated
by the right invariant vector fields $R_j$ on $SO(3)$ such that
\[
[L_i, R_j]=0, \quad [R_i, R_j]=-\frac{1}{2}\epsilon_{ijk}R_k.
\]
These vector fields can be expanded in terms of $L_i$s and therefore
also in terms of the dual basis $E_1, E_2, E_3$ as
\[
R_i={M^j}_i E_j
\]
where the components of the matric $M$ are explicit functions of
$(r, \psi, \theta, \phi)$, and we have verified  that
${\mathcal{L}}_{R_i} \Phi=0$.
Combining this with (\ref{pauls_c}) leads
to three first integrals
\[
p_i\equiv {R_i}^a(u_a+\kappa\Phi_a)= {M^j}_i (u_j+\kappa\Phi_j).
\]
These are three linear equations for ${\bf u}$ which can be soved 
as ${\bf u}=M^{-1}{\bf p}-\kappa{\bf \Phi}$, or (using the explicit form of $M$)
\begin{eqnarray}
\label{uspaul}
u^1&=&\frac{\kappa}{2}c_1 rf+\frac{1}{r}\Big(p_2(\cos{\phi}\cos{\psi}-
\sin{\phi}\sin{\psi}\cos{\theta})-p_3
(\sin{\phi}\cos{\psi}-
\cos{\phi}\sin{\psi}\cos{\theta})+p_1\sin{\theta}\sin{\psi}\Big)\nonumber\\
u^2&=&\frac{\kappa}{2}c_2 rf+\frac{1}{r}\Big(-p_2(\cos{\phi}\sin{\psi}-
\sin{\phi}\cos{\psi}\cos{\theta})+p_3
(\sin{\phi}\sin{\psi}-
\cos{\phi}\cos{\psi}\cos{\theta}) +p_1\sin{\theta}\cos{\psi} \Big)\nonumber\\
u^3&=&\frac{\kappa}{2f}c_3 r+\frac{1}{rf}(p_1\cos{\theta}+
p_2\sin{\phi}\sin{\theta}+p_3\cos{\phi}\sin{\theta}).
\end{eqnarray}
\subsection{A final system of 1st order ODEs}
Comparing the expressions  (\ref{uspaul})
with (\ref{uscor}) gives a system
of three 1st order equations for three unknown functions
$\theta(s), \psi(s), \phi(s)$
\begin{subequations}
\begin{eqnarray}
\label{finaleq}
\dot{\phi}&=& \frac{1}{r^2}(2p_1-2\cot{\theta}(p_2\sin{\phi}+p_3\cos{\phi}))
+\kappa f\frac{c_1\sin{\psi}+c_2\cos{\psi}}{\sin{\theta}},\\
\dot{\theta}&=&\frac{1}{r^2}(p_2\cos{\phi}-p_3\sin{\phi})+\kappa f (c_1\cos{\psi}
-c_2\sin{\psi}),
 \\
\dot{\psi}&=&\kappa\Big(f\cot{\theta}(c_1\sin{\psi}+c_2\cos{\psi})
+\frac{1}{f^2}c_3)\Big)+\frac{2}{r^2f^2\sin{\theta}}(p_2\sin{\phi}+p_3\cos{\phi})\\
&&+\frac{2\alpha^4}{r^6f^2\sin{\theta}}(\cos{\theta}^2(p_2\sin{\phi}+p_3\cos{\phi})
-p_1\cos\theta\sin{\theta}).\nonumber
\end{eqnarray}
\end{subequations}
The equation for $\dot{r}$ can now be obtained from the 1st integral $g(u,  u)=1$.
\subsection{An integrable case:  circles on $SO(3)$ orbits}
In this section we shall establish the complete integrability of the
of the conformal geodesic equations under additional assumption
that the conformal geodesics lie on surfaces of constant $r$.

\noindent
{\bf Proof of Proposition \ref{propEH}.}
The circle equations can be reduced to a quadrature under the additional 
assumption that $\dot{r}=0$, that is the circles are confined to
the orbits of the isometric $SO(3)$ action which are surfaces of constant $r$.
The condition $\dot{r}=0$ is equivalent to $u^4=0$.  Set
\[
u^1=\alpha(s), \quad u^2=\beta(s), \quad u^3=\gamma(s),\quad
\frac{k^4}{r\sqrt{1-k^4}}=R=\mbox{const}. 
\]
Equations (\ref{ueq1}) take the form
\begin{subequations}
\label{ueq2}
\begin{eqnarray}
\dot{\alpha}&=&R\;\beta\gamma+c_2\gamma-c_3\beta\\
\dot{\beta}&=& -R\;\alpha\gamma+c_3\alpha-c_1\gamma\\
\dot{\gamma}&=&c_1\beta-c_2\alpha.
\end{eqnarray}
\end{subequations}
Algebraic manipulations with these equations give expressions
for $\alpha$ and $\beta$ in terms of $\gamma$ and its derivatives
\begin{eqnarray}
\label{alphabeta}
\alpha&=&-\frac{c_1}{({c_1}^2+{c_2}^2)(R\gamma-c_3)}
\ddot{\gamma}
-\frac{c_2}{{c_1}^2+{c_2}^2}\dot{\gamma}
-\frac{c_1}{R\gamma-c_3} \gamma\nonumber\\
\beta&=&-\frac{c_2}{({c_1}^2+{c_2}^2)(R\gamma-c_3)}
\ddot{\gamma}
+\frac{c_1}{{c_1}^2+{c_2}^2}\dot{\gamma}
-\frac{c_2}{R\gamma-c_3} \gamma.
\end{eqnarray}
The last equation in (\ref{ueq2}) now holds identically. The first two equations
reduce to a single 3rd  ODE 
\be
h\;\dddot{h}=(\ddot{h}-h^3+m_0)\dot{h}, \quad\mbox{where}
\;\; h(s)=R\gamma(s)-c_3, \quad\mbox{and}\;\; m_0=c_3({c_1}^2+{c_2}^2). 
\ee
This equation can be integrated in terms of elliptic functions
\be
\label{diff2ce}
s-m_3=2\int\frac{dh}{\sqrt{-h^4+4m_1h^2-8m_0h+4m_2}},
\ee
where $m_0, m_1, m_2, m_3$ are constants of integration. Inverting this
formula for $h$, and substituting in the expressions above gives
$\alpha(s), \beta(s)$. The elliptic functions are however
not neccesary to describe the conformal circles. To see this differentiate
(\ref{diff2ce}) twice, to obtain expressions for $\dot{h}$ and $\ddot{h}$
in terms of $h$. These expressions are then substituted into 
(\ref{alphabeta}), and (remembering
that $\gamma=R^{-1}(h(s)+c_3)$) yield the explicit parametrisation
${\bf u}(h)$
of $u^1, u^2, u^3$ in terms of a parameter $h$, and five integration constant
(note that there is one constraint on the constants
$m_0, m_1, m_2, m_3, c_1, c_2$ coming from $g(u, u)=1$). The formulae
for  ${\bf u}(h)$ are now substituted into (\ref{uspaul}) which yields
relations of the form
\[
G_1(\phi, \psi, \theta, h)=0, \quad G_2(\phi, \psi, \theta, h)=0,\quad
G_3(\phi, \psi, \theta, h)=0
\]
where the functions $G_i$ can be read--off from (\ref{uspaul}), and
depend on the 8 constants of integrations (five from the procedure above,
and $p_1, p_2, p_3$). The parameter $h$ may be elliminated between these three
relations which gives unparametrised conformal circles
on the $r=$const orbits in the Eguchi--Hanson space. They are given by
(\ref{starstar}).
\koniec
\subsection{Separability of the Hamilton--Jacobi equation}
In this section we shall find another integrable subcase of the conformal
geodesic equations on the Eguchi--Hanson background. This will be done by
separating the associated Hamilton--Jacobi equation.
Following \cite{gibbons_ruback} again, we introduce prolate spheroidal coordinates
adapted to the potential (\ref{VEH})
\[\zeta=\frac{1}{2\alpha}(r_1+r_2),\;\;\;\lambda=\frac{1}{2\alpha}(r_1-r_2),\]
so these are constant on confocal ellipses and confocal hyperbolae respectively. Introduce the angles $\theta_1,\theta_2$ as the angles supported at $p_1,p_2$ respectively w.r.t. the $z$-axis; now
\[r_1=\alpha(\zeta+\lambda),\;\;r_2=\alpha(\zeta-\lambda),\;\;\cos\theta_1=\frac{\zeta\lambda+1}{\zeta+\lambda},\;\;\cos\theta_2=\frac{\zeta\lambda-1}{\zeta-\lambda}.\]
We need $V$ and $\omega$ in these coordinates: for the Eguchi-Hanson metric
\[V=\frac{1}{r_1}+\frac{1}{r_2}=\frac{2\zeta}{\zeta^2-\lambda^2},\;\;\omega=\Omega d\phi\quad\mbox{with}\quad\Omega=-2\frac{\lambda(\zeta^2-1)}{{\zeta^2-\lambda^2}},\]
with conventions
\[\alpha(\zeta^2-1)V_\zeta=\Omega_\lambda,\;\;-\alpha(1-\lambda^2)V_\lambda=\Omega_\zeta.\]
Claim
\[d\rho^2+dz^2+\rho^2d\phi^2=\alpha^2(\zeta^2-\lambda^2)\left(\frac{d\zeta^2}{\zeta^2-1}+\frac{d\lambda^2}{1-\lambda^2}\right)+\alpha^2(\zeta^2-1)(1-\lambda^2)d\phi^2\]
whence
\[g_{EH}=V\left(\alpha^2(\zeta^2-\lambda^2)\left(\frac{d\zeta^2}{\zeta^2-1}+\frac{d\lambda^2}{1-\lambda^2}\right)+\alpha^2(\zeta^2-1)(1-\lambda^2)d\phi^2\right)+V^{-1}(d\psi+\Omega d\phi)^2.\]
From the geodesic equations we obtain constants of the motion
\[p_\psi=V^{-1}(\dot\psi+\Omega\dot\phi):=E,\;\;p_\phi={\alpha^2}V(\zeta^2-1)(1-\lambda^2)\dot\phi+\Omega E:=J,\]
from the Killing vectors $K:=\partial_\psi$ and $L:=\partial_\phi$, so that
\[p_\phi-\Omega p_\psi=\alpha^2V(\zeta^2-1)(1-\lambda^2)\dot\phi.\]
Seek a solution
\[S=E\psi+J\phi+F(\zeta)+G(\lambda)\]
for the Hamilton-Jacobi equation then
\[\mu^2=\frac{\zeta^2-1}{{\alpha^2}V(\zeta^2-\lambda^2)}F_\zeta^2+\frac{1-\lambda^2}{{\alpha^2}V(\zeta^2-\lambda^2)}G_\lambda^2+\frac{(J-\Omega E)^2}{{\alpha^2}V(\zeta^2-1)(1-\lambda^2)}+VE^2.\]
Multiply by ${\alpha^2}V(\zeta^2-\lambda^2)$ and simplify to obtain
\be\label{hj1}2\alpha\zeta=(\zeta^2-1)F_\zeta^2+(1-\lambda^2)G_\lambda^2+J^2(\frac{1}{\zeta^2-1}+\frac{1}{1-\lambda^2})
+\frac{4JE\lambda}{1-\lambda^2}+\frac{4E^2}{1-\lambda^2},\ee
which separates.

Up to this point we are following \cite{gibbons_ruback}. To extend to the Lorentz force law, we need a Maxwell potential for $F_{ab}$ which commutes with $K$ and $L$. Now there is a problem: for anti-self-dual Taub-NUT there was no loss in generality in choosing $F=-\Omega^3$ but in Eguchi-Hanson there is, as $K$ will Lie-drag any constant $\phi$ but $L$ will not. Thus we can obtain integrability for a sub-class of conformal geodesics, those for which $F=-\Omega^3$, but not for all. With this restriction then we have
\[\Omega^3=-(d\psi+\omega)\wedge dz+Vdx\wedge dy.\]
If the potential for $-\Omega_3$ is
\[\Phi=-z(d\psi+\omega)+Xd\phi=-z(d\psi+\Omega d\phi)+Xd\phi,\]
then
\[dX\wedge d\phi=zd\Omega\wedge d\phi-V\rho d\rho\wedge d\phi.\]
This can be solved by
\[X=z\Omega+2\alpha\zeta\]
when
\[\Phi=(X-z\Omega)d\phi-zd\psi=2\alpha\zeta d\phi-{\alpha}\zeta\lambda d\psi.\]
To modify the Hamilton-Jacobi equation we make the replacements
\[J\rightarrow J+2e\alpha\zeta,\;\;E\rightarrow E-e\alpha\zeta\lambda\]
into (\ref{hj1}). The added terms are
\[4e\alpha J\frac{\zeta^3}{\zeta^2-1}+4e^2\alpha^2\frac{\zeta^4}{\zeta^2-1},\]
so this still separates, and this sub-class of the conformal geodesic equations is completely integrable.

\section{Circles on $\CP^2$}
\label{section3}
Consider the Fubini--Study metric on $\CP^2$. It is Einstein with the Ricci 
scalar equal to $24$, has anti--self--dual (ASD) Weyl tensor, and is K\"ahler,
but with the ASD K\"ahler two-form (this is some times referred to 
as `opposite orientation'). The local form of the metric is (see e.g. \cite{GP})
\be
\label{cp2metric}
g_{\CP^2}=\frac{dr^2}{(1+r^2)^2}+
\frac{1}{4}\frac{r^2\sigma_3^2}{(1+r^2)^2}
+\frac{1}{4}\frac{r^2}{1+r^2}(\sigma_1^2+\sigma_2^2).
\ee
The metric is regular everywhere on $\CP^2$, and 
the apparent singularity at $r=0$ results from 
using spherical polar coordinates. The metric is also regular at $r=\infty$. To see this set 
$r=\rho^{-1}$. Fixing $(\phi, \theta)$ now gives $g\sim d\rho^2+(1/4)\rho^2 d\psi^2$
near $\rho=0$. This is a removable bolt singularity (see \cite{Dbook}) if $\psi$ is periodic with the period $2\pi$. At $\rho=0$ the three--dimensional orbits of $SU(2)$ collapse to a two--sphere of constant radius. We shall refer to this as the $\CP^1$ at
infinity.

It was shown in \cite{DT} that in case of 
ASD Einstein manifolds 
with non--zero Ricci scalar there is a one-to-one correspondence between Killing vectors, and self--dual CKY tensors: If $K$ is a Killing vector, then 
its self dual derivative 
\[
Y\equiv \frac{1}{2}(dK+*dK)
\]
satisfies (\ref{CKY}). The Fubini--Study metric is a symmetric space
$M=SU(3)/U(2)$, so it has 8 Killing vectors generating the Lie algebra
$\mathfrak{su}(3)$, and therefore 8 self--dual CKYs. This gives rise to 8 first
integrals of the form (\ref{Q}). The 9th first integral is given by the parallel
CKY which is the anti--self--dual K\"ahler form. In Section  we aim to use 
these first integrals to find all conformal geodesics on $\CP^2$.
We shall first prove
\begin{prop}
\label{prop_CP2}
All conformal geodesics on the Fubini--Study metric (\ref{cp2metric})
on $\CP^2$ on the surfaces of constant $r$ are trajectories of the Killing
vector $\p/\p\psi$, or are of the form
\be
\label{theta_phi_psi_new}
\chi={\rm arccot}
\Big(\frac{\kappa(c_2\cos{(\kappa s)}-\sin{(\kappa s))}}{Q(c_3+c_2\sin{(\kappa s)}+\cos{(\kappa s))}}\Big), \quad
\theta={\rm arccot}\Big({\frac{Q-\dot{\chi}}{P\sin{\chi}}}\Big), \quad \phi=P\int\frac{\sin{\chi}}{\sin{\theta}}ds
\ee
where
\[
Q=2\gamma/r-\gamma r,\quad  P=\frac{2\sqrt{(1+r^2)(1-\gamma^2)}}{r}, \quad \kappa=\sqrt{P^2+Q^2}
\]
and $r, \gamma, c_1, c_2, c_3$ are constants.
\end{prop}\noindent
and then deduce the the general form of conformal geodesics:
\begin{theo}
\label{theo_CP2}
All conformal geodesics of the Fubini--Study metric on $\CP^2$
are of the form ${\bf \rho}\circ\Gamma$ where $\Gamma$ is some conformal 
geodesic on the $r$=const surface given in Proposition \ref{prop_CP2}, 
and ${\bf \rho}$ is an isometry of the Fubini--Study metric.
\end{theo}
\noindent
{\bf Proof of Proposition  \ref{prop_CP2}}.
Pick a tetrad of one--forms
\[
e^1=\frac{r}{2\sqrt{1+r^2}}{\sigma_1}, \quad
e^2=\frac{r}{2\sqrt{1+r^2}}{\sigma_2}, \quad
e^3=\frac{r}{2(1+r^2)}\sigma_3, 
\quad 
e^4=\frac{dr}{1+r^2},
\]
so that the metric (\ref{cp2metric}) is $g=\delta_{ab}e^ae^b$,
and set 
\be
\label{formulaua}
u={\alpha} E_1+{\beta} E_2+{\gamma} E_3+{\delta} E_4, \quad 
a=a^1 E_1+a^2 E_2+a^3 E_3+a^4 E_4
\ee
where $E_a, a=1, 2, 3, 4$ is the dual tetrad of vector fields, and 
$|u|^2=1$.
The conformal geodesic equations (\ref{conf_circ_e}) become
\begin{subequations}
\label{fs}
\begin{eqnarray}
a^1&=&\dot{\alpha}+2 r {\beta}{\gamma}+\frac{1}{r}{\alpha} {\delta}
\label{fs1}\\
a^2&=&\dot{\beta}-2 r {\gamma} {\alpha} +  \frac{1}{r}{\beta} {\delta}\label{fs2}\\
a^3&=&\dot{\gamma}-r {\gamma}{\delta}+\frac{1}{r}{\gamma}{\delta}\label{fs3}\\
a^4&=&\dot{\delta}+r {\gamma}^2-\frac{1}{r}(1-{\delta}^2)\label{fs4}
\end{eqnarray}
\end{subequations}
and
\begin{subequations}
\begin{eqnarray}
\dot{a}^1&=&\frac{1}{r}({\beta}a^3-{\alpha}a^4-{\gamma}a^2)-2r {\gamma}a^2-|a|^2 {\alpha}\\
\dot{a}^2&=&\frac{1}{r}({\gamma}a^1-{\beta}a^4-{\alpha}a^3)+2r {\gamma}a^1-|a|^2 {\beta}\\
\dot{a}^3&=&\frac{1}{r}({\alpha}a^2-{\gamma}a^4-{\beta}a^1)+r {\gamma}a^4-|a|^2 {\gamma}\\
\dot{a}^4&=&\frac{1}{r}({\alpha}a^1+{\beta}a^2+{\gamma}a^3)-r {\gamma}a^3-|a|^2 {\delta},
\end{eqnarray}
\end{subequations}
together with
\[
|u|^2=1, \quad \delta_{ef}a^ea^f=|a|^2, \quad\mbox{where}\;\;u=(\alpha, \beta, \gamma, \delta).
\]
Note that $g(u, a)=0$ is satisfied identically if $u^a$ is  unit.

The CKY corresponding to the Killing vector $\p/\p \psi$ 
is
\[
\frac{r^2}{r^2+1}\Big(e^1\wedge e^2+e^3\wedge e^4\Big),
\]
and gives rise to a first integral
\be
\label{CYint}
C_Y:=\frac{r^2}{r^2+1}\Big({\alpha}a^2 -{\beta}a^1 +{\gamma}a^4-{\delta}a^3\Big)+\frac{2r}{r^2+1}{\gamma}.
\ee
The ASD K\"ahler form is
\be
\label{kahler_form}
e^1\wedge e^2-e^3\wedge e^4.
\ee
This gives rise to a first integral
\be
\label{CKint}
C_K:={\alpha}a^2 -{\beta}a^1 -{\gamma}a^4+{\delta}a^3.
\ee
Combining the two gives
\begin{subequations}
\begin{eqnarray}
(r^2+1)C_Y-r^2C_K&=&2r^3{\gamma}({\gamma}^2+{\delta}^2)+ 
2r^2({\gamma}\dot{\delta} -{\delta}\dot{\gamma}). \label{1stint1}\\
(r^2+1)C_Y+r^2C_K&=&4r^3{\gamma}({\gamma}^2+{\delta}^2)+ 
2r^2({\alpha}\dot{\beta} -{\beta}\dot{\alpha})
+2r{{\gamma}}-4r^3{\gamma}.\label{1stint2}
\end{eqnarray}
\end{subequations}
The first integral $|a|^2$ after some algebra simplifies to
\begin{eqnarray}
\label{AAeq}
|a|^2&=&\delta_{ef}\dot{u}^e\dot{u}^f+r^2{\gamma}^2(4-3({\gamma}^2+{\delta}^2))
+2r{\gamma}(2{\beta}\dot{\alpha}-2{\alpha}\dot{\beta}+
{\gamma}\dot{\delta}-{\delta}\dot{\gamma})\nonumber\\
&& -2{\gamma}^2 -\frac{2}{r}\dot{\delta}+\frac{1}{r^2}(1-{\delta}^2).
\end{eqnarray}
\vskip 5pt
We  now aim to find all circles such that ${\delta}=0$, or equivalently all circles lying on surfaces of constant $r$.
The vector field $u$ is unit, so set $\beta=\sqrt{1-\alpha^2-\gamma^2}$. Equations (\ref{1stint1}, \ref{1stint2}) and (\ref{AAeq}) give
\begin{subequations}
\begin{eqnarray}
\label{integralscp}
0&=&2(r\gamma)^3-r^2(C_Y-C_K)-C_Y \label{integralscpa} \\
0&=&2r^2(\alpha\dot{\beta}-\beta\dot{\alpha})+4(r\gamma)^3-4r^3\gamma+2r\gamma
-r^2(C_K+C_Y)-C_Y \label{integralscpb}\\
A&=&\dot{\alpha}^2+\dot{\beta}^2+(r\gamma)^2(4-3\gamma^2)+4r\gamma(\beta\dot{\alpha}-\alpha\dot{\beta})-2\gamma^2+\frac{1}{r^2}\label{integralscpc}
\end{eqnarray}
\end{subequations}
so that $\gamma$ is a constant along circles (and this constant depends on $r$).
We also get $a^4=r\gamma^2-1/r$, which is a constant along circles, 
and $a^3=0$. 
Substituting $\beta=\sqrt{1-\alpha^2-\gamma^2}$ into 
(\ref{integralscpb}) yields a 1st order equation for $\alpha(s)$ which we can integrate to find
\begin{eqnarray}
\label{cp2cir}
\alpha(s)&=&\sqrt{1-\gamma^2}\sin{(Bs+c_1)}, \quad 
\beta(s)=\sqrt{1-\gamma^2}\cos{(Bs+c_1)}, \quad\mbox{where}\\
B&=&
\frac{4r^3(\gamma^3-\gamma)-r^2(C_K+C_Y)+2 r\gamma-C_Y}{2r^2(1-\gamma^2)},
\quad\mbox{where}\;\;\gamma\neq\pm 1\nonumber
\end{eqnarray}
and $\gamma$ is a root of 
(\ref{integralscpa}). 

To find any additional constraints on $B$ we substitute $u=(\alpha(s), \beta(s), \gamma, 0)$ into the conformal geodesic equations. We find that these equations
hold iff 
\be
\label{conditons_for_B}
B=-3\gamma r, \quad\mbox{or}\quad \gamma^2r^2+\gamma B+1=0.
\ee
The squared norm of acceleration (\ref{integralscpc})
is constant along conformal circles, and is
given by
\[
|a|^2=B^2(1-\gamma^2)+4r\gamma B(1-\gamma^2)+r^2\gamma^2(4-3\gamma^2)-2\gamma^2+
\frac{1}{r^2}.
\]
We can now compute the {\em complex torsion} $\tau=C_K/|a|$ and find that
it can take any value between $(-1, 1)$ if $B=-3\gamma r$, but $\tau^2=1$
if $B=-\gamma r^2-\gamma^{-1}$. We shall therefore take $B=-3\gamma r$ from now on, and return to this point in the proof of Theorem \ref{theo_CP2}.

If $\gamma=\pm 1$, then the circle is a trajectory of
the Killing vector $\p/\p \psi=r/(r^2+1)E_3$, and $\alpha=\beta=0$
with the  corresponding acceleration  $a=(r^2-1)/r E_4$. 
In this case we  find that $u=E_3$ and 
\be
\label{psigeod}
\nabla_u \nabla_u u=-\frac{(r^2-1)^2}{r^2}u,
\ee
so the circle equations (\ref{conf_circ_e}) hold\footnote{We also note that in general $\dot{r}=(1+r^2){\delta}$, and use
this to verify that $u=E_a (a=1, 2, 3, 4)$ satisfies the circle eq. In the last case we get $|a|^2=0$.} with $|a|^2=(r^2-1)^2/r^2$.

\vskip5pt
To complete the calculation for $\gamma\neq \pm 1$ we need to find the dependence of
the coordinates $(\phi, \theta, \psi, r)$ on the arc-lengh parameter $s$ along 
the circles (\ref{cp2cir}). Setting
\[
u=\dot{\phi}\frac{\p}{\p\phi}+\dot{\theta}\frac{\p}{\p\theta}+
\dot{\psi}\frac{\p}{\p\psi}+\dot{r}\frac{\p}{\p r},
\]
equating it to  $u$ in (\ref{formulaua}) and substituting 
(\ref{cp2cir}) yields a coupled system of three 1st order ODEs
\begin{subequations}
\begin{eqnarray}
\label{system3}
\dot{\phi}&=&2\frac{\sqrt{1+r^2}}{r}\frac{\sqrt{1-\gamma^2}}{\sin{\theta}}
\sin{(\psi+Bs+c_1)},\label{system3a} \\
\dot{\theta}&=&2\frac{\sqrt{1+r^2}}{r}{\sqrt{1-\gamma^2}}\cos{(\psi+Bs+c_1)},\label{system3b}\\
\dot{\psi}&=&2\frac{1+r^2}{r}\gamma
-2\frac{\sqrt{1+r^2}}{r}\sqrt{1-\gamma^2}\cot{\theta}
\sin{(\psi+Bs+c_1)}\label{system3c},\\
\dot{r}&=&0\label{system3d} .
\end{eqnarray}
\end{subequations}
To simplify this slightly, set
\be
\label{introchi}
\chi=\psi+Bs+c_1, \quad P=\frac{2\sqrt{(1+r^2)(1-\gamma^2)}}{r}=\mbox{const}, 
\quad
Q=B+2\frac{1+r^2}{r}\gamma=\mbox{const}
\ee
so that
\be
\label{simpsimp}
\dot{\phi}=P\frac{\sin{\chi}}{\sin{\theta}}, \quad
\dot{\theta}=P\cos{\chi}, \quad \dot{\chi}=Q-P\cot{\theta}\sin{\chi}.
\ee
We claim that this system is solvable by quadratures. Given $\chi(s)$ we find
\be
\label{theta_phi}
\theta={\rm arccot}\Big({\frac{Q-\dot{\chi}}{P\sin{\chi}}}\Big), \quad
\mbox{and then}\quad \phi=P\int\frac{\sin{\chi}}{\sin{\theta}}ds.
\ee
The final equation for $\chi$ arises from the middle equation
in (\ref{simpsimp}):
\[
 P^2\cos{(\chi)}^{3}
-2\cos{(\chi)}\dot{\chi}^2+
3Q \cos{(\chi)}\dot{\chi}-
\cos{(\chi)}({P}^{2}+Q^2)+
\sin{(\chi)} \ddot{\chi}=0.
\]
The general solution of this equation is given by
\be
\label{formchi}
\chi={\rm arccot}
\Big(\frac{\kappa(c_2\cos{(\kappa s)}-\sin{(\kappa s))}}{Q(c_3+c_2\sin{(\kappa s)}+\cos{(\kappa s))}}\Big), 
\quad \mbox{where}\quad \kappa= \sqrt{P^2+Q^2}.
\ee
This can now be substituted in (\ref{theta_phi}) to find the explicit form of 
$\theta$. It turns out that (according to MAPLE) even the $\phi$ integral
can be computed in terms of elementary functions, but the answer is not 
illuminating. The integration will in any case introduce another constant - 
call it $c_4$. Thus there exists a seven--dimensional
family of conformal circles on surfaces of constant $r$. This family is
parametrised by $c_1, c_2, c_3, c_4, C_Y, C_K, r$, and explicitly given 
by formulae (\ref{formchi}, \ref{theta_phi}, \ref{introchi}), 
where $B$ is given by (\ref{cp2cir}), and $\gamma$ is given by (\ref{integralscpa}).
\koniec
We shall now use the conformal geodesics on the $r=$ const surfaces
and conjugate them with the elements of the isometry group $SU(3)$ to find
all conformal geodesics on $\CP^2$. The fact that all conformal geodesics
arise from this construction will follow from the results of \cite{MO} which we 
now recall.

We shall say that two circles $\Gamma_1$ and $\Gamma_2$ are congruent  if there exists
a holomorphic isometry (note - all isometries on $\CP^2$ are holomorphic) $
{\bf \rho}$ such that $\Gamma_1={\bf \rho}\circ \Gamma_2$.

In  \cite{MO} (Theorem 5.1) it has been shown that 
two circles on $\CP^2$ with the same complex torsion 
\be
\label{torsion}
\tau=g_{\CP^2}(J(u), a/|a|)
\ee
and curvature $|a|^2$ are congruent. Here $J$ is the complex structure
of the ASD K\"ahler form on $\CP^2$ and we note that $\tau=C_K/|a|$, where
$C_K$ is the first integral (\ref{CKint}). We also note 
$\tau^2\leq |J(u)|^2 |a|^2/|a|^2=1$,
so that $\tau \in [-1, 1]$.
Doing the counting, at a point on 
$r$=constant there are $2$
dimensions of (unit) velocity $u$ and then $2$ dimensions of acceleration orthogonal to velocity so with $+3$ for the dimension of $r$=constant surface and $-1$ for the one-dimension of the
curve that is $6$ dimensions of conformal geodesics confined to $r$=const. But there is a four--dimensional space of Killing vectors tangent to $r$=constant to move them about, leaving $2$
parameters. In the proof of Theorem \ref{theo_CP2} we shall show that  these are $|a|$ and $\tau$. 

The general counting for unparametrised
conformal geodesics on $\CP^2$ is as follows:
There are 
$3$ degrees of freedom of velocity, $3$ of acceleration, $4$ of position
minus $1$ for the dimension of the curve. This gives
gives a $9$ dimensional space of conformal geodesics on $\CP^2$. 
There are $8$
Killing vectors, but an orbit of $SU(3)$ has fixed $|a|$ and 
$\tau$ which is $7$ dimensions in the space of conformal
geodesics. Therefore there is always $1$-dimensional isometry group 
stabilising a conformal geodesic i.e. every conformal geodesic 
$\Gamma(s)=\mbox{exp}(sK) \;\Gamma(0)$,
is a 
trajectory of  some Killing vector $K$
which confirms the results of \cite{maeda_kil}. 
This counting is correct as long as the Killing vector $K$ is tangent to the conformal
geodesic, rather than being vanishing on it. To  examine this we look at fixed points of Killing vectors in $\CP^2$. A one parameter subgroup of the isometry group comes from an $SU(3)$ matrix,
and fixed points correspond to eigenvectors; if the eigenvalues are distinct
then the Killing vector vanishes at isolated points (like 
$\p_\phi$ in the coordinate system (\ref{cp2metric})) but if there is  a
repeat then the Killing vector vanishes on a $\CP^1$ in $\CP^2$ (like $\p_\psi$ which vanishes
at $r=\infty$ - the $\CP^1$ at infinity). So if a conformal geodesic  does not lie on a $\CP^1$, then it
can not lie in the zero set of a Killing so must be a Killing vector trajectory. If it does lie
on a $\CP^1$ then there is a Killing vector that vanishes everywhere along the geodesic but it does not matter; take it to be the $\CP^1$ 
at infinity which is a metric sphere and
extrinsically flat, so that the conformal geodesics are genuine
circles and are Killing vector trajectories.
\vskip5pt
\noindent
{\bf Proof of Theorem \ref{theo_CP2}.}
In view of the discussion above, and the results of \cite{MO} it is sufficient
to show that geodesics on constant $r$ surfaces
from Proposition \ref{prop_CP2} can have arbitrary values
of the torsion $\tau$, and the curvature $|a|^2$. We shall first
consider the case when $\tau=\pm 1$. Then the trajectory of the Killing
vector $\p/\p\psi$ is a conformal geodesic (\ref{psigeod})
with $|a|^2=(r^2-1)^2/r^2$ which can take any value in $[0, \infty)$.

We shall therefore assume that $\tau^2<1$, and show that for any value
of $(\tau, |a|)$ there is a corresponding value of $(\gamma, r)$.
This will be achieved  by finding two relations
between $(\tau, |a|)$ and  $(\gamma, r)$.

Recall from (\ref{conditons_for_B}) that
$B=-3\gamma r$, or $\gamma^2r^2+\gamma B+1=0$. In the latter case
we find ${C_K}^2=|a|^2$, so that $\tau=\pm 1$, but we already
have a conformal geodesics with these values of $\tau$. We shall therefore
assume that
\be
\label{bminusthree}
B=-3\gamma r
\ee
in which case we find
\be
\label{p4}
|a|^2=\gamma^2 r^2-2\gamma^2+\frac{1}{r^2}
\ee
and
\be\label{p5}
C_K=\frac{\gamma}{r}(1+r^2-2\gamma^2r^2).\ee
As a check, calculate
\[1-\tau^2=1-\frac{C_K^2}{|a|^2}=\frac{(1-\gamma^2)(1-2r^2\gamma^2)^2}{(1-\gamma^2)+\gamma^2(r^2-1)^2}>0,\]
as expected. Also (\ref{p5}) will give the sign of $C_K/\gamma$ (note that, given a conformal geodesic, then traversing it in the opposite direction switches the signs on $C_K$ and $\gamma$). Solving (\ref{p4}) for $\gamma^2$
gives 
\be\label{p6}\gamma^2=\frac{r^2|a|^2-1}{r^2(r^2-2)},\ee
(once we have $r^2$ we will  need to check that this is less than 1.)
This relation together with 
\[
C_K=-r^3\gamma^3+r\gamma(1+|a|^2)
\]
(which follows from (\ref{CKint})) gives a cubic
\[
F(X)\equiv(XA-1)(X-(2A+1))^2-TA(X-2)^3=0, \quad
\mbox{where}\quad  X\equiv r^2, A\equiv |a|^2, T\equiv \tau^2.
\]
We will argue that this cubic has 
three positive roots which coalesce at $X=2$ when $A=1/2$.

\vskip5pt

Let us deal separately with $r^2=2$: from (\ref{p4}) for this case $|a|^2=1/r^2=1/2$ for any $\gamma^2$; then from (\ref{p5})
 \[\tau=\frac{C_K}{|a|}=\gamma(3-4\gamma^2),\]
and for any $\tau$ with $-1<\tau<1$, which is the currently allowed range, inspection of this cubic shows that there are three roots for $\gamma$ in $(-1,1)$. 

Now we can suppose that $r^2\neq 2$. 
Substitute for $\gamma^2$ into the square of (\ref{p5}):
 \[C_K^2=\frac{\gamma^2}{r^2}(1+r^2-2\gamma^2r^2)^2=\frac{r^2|a|^2-1}{r^4(r^2-2)}\left(1+r^2-2r^2\left(\frac{r^2|a|^2-1}{r^2(r^2-2)}\right)\right)^2 \]
 \[=\frac{(r^2|a|^2-1)}{(r^2-2)^3}(r^2 -(1+2|a|^2))^2. \]
 Rationalise:
 \be\label{eq1}
 0=(r^2|a|^2-1)(r^2 -(1+2|a|^2))^2-C_K^2(r^2-2)^3,\ee
 and put $r^2=X$ to define the cubic
 \begin{eqnarray}\label{eq2}F(X)&:=&X^3(|a|^2-C_K^2)-X^2(1+2|a|^2+4|a|^4-6C_K^2)\\
 &&+X(2+5|a|^2+4|a|^4+4|a|^6-12C_K^2)\nonumber\\
 &&-(1+4|a|^2+4|a|^4-8C_K^2),\nonumber
 \end{eqnarray}
 whose vanishing we want. Note that, in $F$:
 \begin{itemize}
  \item The coefficient of $X^3$ is
 \[|a|^2-C_K^2=(1-\tau^2)|a|^2,\]
 which is strictly positive;
 \item 
 the coefficient of $X^2$ is
 \[-(1+2|a|^2(1-3\tau^2)+4|a|^4)=-((1-2|a|^2)^2+6|a|^2(1-\tau^2))\]
 which is strictly negative; 
 \item the coefficient of $X$ is 
 \[2+5|a|^2+4|a|^4+4|a|^6-12C_K^2>2-7|a|^2+4|a|^4+4|a|^6\]\[=(2|a|^2-1)^2(|a|+2),\]
 where the first inequality uses $\tau^2<1$, so that this coefficient is strictly positive;
\item and the constant term is
 \[-(1+4|a|^2+4|a|^4-8C_K^2)=-(1-2|a|^2)^2-8|a|^2(1-\tau^2),\]
which is strictly negative.
 \end{itemize}
 Thus $F(X)$ has no negative roots and at least one positive root. We will see below there are three positive roots, which coalesce at $r^2=2$ when $|a|^2=1/2$.

Next we note some particular values of $F$:
\be\label{eq3}F(1/|a|^2)=\frac{\tau^2}{|a|^8}(2|a|^2-1)^3,\;\;F(2)=(2|a|^2-1)^3\ee
and
\be\label{eq4}F(1+2|a|^2)=-C_K^2(2|a|^2-1)^3.\ee
Now we can analyse the roots of $F$
\begin{itemize}
\item If $|a|^2=1/2$ then $r^2=2$ is a root, and by inspection of (\ref{eq1}) we see that it is three times repeated. By the discussion above, there is now a value of $\gamma$ giving any value of $\tau$.
\item If $2|a|^2-1<0$ then the evaluations in (\ref{eq3}) are both negative and that in (\ref{eq4}) is positive, while
\[1+2|a|^2<2<1/|a|^2.\]
Therefore there is one root in each of the ranges $(0,1+2|a|^2),(1+2|a|^2,2)$ and $(1/|a|^2,\infty)$. All three make $\gamma^2$ in (\ref{p6}) positive but we shall choose the largest root $r_3^2$ which has
\[r_3^2>1/|a|^2>2\]
so that there are positive $\delta,\epsilon$ with $\delta<\epsilon$ and
\[r_3^2=1/|a|^2+\delta=2+\epsilon.\]
Now from (\ref{p6})
\[\gamma^2=\frac{r_3^2|a|^2-1}{r_3^2(r_3^2-2)}=\frac{\delta|a|^2}{\epsilon r_3^2}<\frac{\delta}{2r_3^2\epsilon}<\frac{\delta}{4\epsilon}<\frac14<1,\]
so that this choice of $r$ leads to a $\gamma$ in the allowed range.
\item If $2|a|^2-1>0$ then the evaluations in (\ref{eq3}) are both positive and that in (\ref{eq4}) is negative, while
\[1/|a|^2<2<1+2|a|^2.\]
Therefore there is one root in each of the ranges $(0,1/|a|^2),(2,1+2|a|^2)$ and $(1+2|a|^2,\infty)$ and all make $\gamma^2$ in (\ref{p6}) positive. Again we choose the largest root $r_3^2$ which therefore has 
\[r_3^2>1+2|a|^2>2.\]
We choose positive $\delta,\epsilon$ with 
\[|a|^2=\frac12+\epsilon,\;\;r_3^2=2+\delta,\]
so that also
\[r_3^2=2+\delta>1+2|a|^2=2+2\epsilon\mbox{  i.e.  }\delta>2\epsilon.\]
Now from (\ref{p6})
\[\gamma^2=\frac{r_3^2|a|^2-1}{r_3^2(r_3^2-2)}=\frac{(1/2+\epsilon)(2+\delta)-1}{\delta(2+\delta)}=\frac14\frac{\delta+4\epsilon+2\epsilon\delta}{\delta+\delta^2/2}     \]
\[=\frac14\left(1+\frac{4\epsilon+2\epsilon^2-(\delta-2\epsilon)^2/2}{\delta+\delta^2/2}\right)<\frac14\left(1+\frac{4\epsilon+2\epsilon^2}{\delta+\delta^2/2}\right)\]
and with $\delta>2\epsilon$:
\[<\frac14\left(1+\frac{4\epsilon+2\epsilon^2}{2\epsilon+2\epsilon^2}\right)<\frac34<1,\]
so again we have $\gamma$ in the allowed range.
\end{itemize}
We have shown that any allowed $(C_K,|a|^2)$ can be obtained from an allowed $(r,\gamma)$.
\koniec
{\bf Remarks.}
\begin{itemize}
\item A smooth curve $\Gamma$ parametrised by an arc--lengh is called a 
helix of proper order $d$, if there exist orthonormal 
vector fields $\{V_1\equiv \dot{\Gamma}, V_2, \dots,  V_d\}$ such that
\[
\nabla_s V_j=-\kappa_{j-1} V_{j-1}+\kappa_j V_{j+1}, \quad\mbox{where}\quad j=1, \dots, d
\]
and $V_0=V_d\equiv 0$, and the functions $\kappa_j=\kappa_j(s)
$ (we set $\kappa_0=\kappa_d\equiv 0$) where $j=1, \dots, d-1$ are constant along $\Gamma$. Thus a helix of order 1 is a geodesic, and a helix of order 2 is a circle. 

Let $\pi: S^{5}\rightarrow \CP^2$
be the Hopf fibration such that the round metric on $S^{5}$ is of the form
\[
g_5=(dt+\Theta)^2+g_{\CP^2},
\]
where $d\Theta=\Omega$ is the K\"ahler form for $\CP^2$ (this is also 
the derivative of $K=\p/\p t$ Killing vector).
In \cite{jap2} it is argued that horizontal (with respect to the $U(1)$
connection $\Theta$) lifts  of circles on $\CP^5$ to $S^{5}$ are
helices of order $2, 3$ or $5$. Moreover the order $2$ occures iff
the first integral (\ref{CKint}) vanishes. 
These are the conformal geodesics of $S^5$ for which velocity and
acceleration are both orthogonal to $K$. A necessary condition is
$J(u,a)=0$ (to see it, differentiate  $g_5(K, a)$ along a conformal geodesic
on $S^5$) which we know can be satisfied, so
the conformal geodesics of $\CP^2$ with $J(u, a)=0$ come from conformal geodesics of $S^5$.
\item There is another special case, where the conformal geodesic equations can be explicitly integrated
with relative ease. Let $J$ be the complex structure of the ASD K\"ahler form
(\ref{kahler_form}), and let $\lambda$ be a non--zero constant. We shall make the ansatz
$a=\lambda J(u)$ so that
\[
a_1=-\lambda\beta, \quad a_2=\lambda \alpha, \quad
a_3=\lambda\delta, \quad a_4=-\lambda\gamma.
\]
Note that the  first integral (\ref{CKint}) is equal to $\lambda$, and
\[
\nabla_u a=\lambda J^2(u)=-\lambda u, \quad \mbox{where}\quad \lambda^2=|a|^2, \quad \mbox{so that}\quad \tau=1.
\]
The circle equations (\ref{conf_circ_e}) hold identically as a consequence
of the definition of $a$, and can be integrated explicitly. We will not reproduce this calculation here,
as we know from Theorem \ref{theo_CP2}, that the resulting conformal geodesics lie
in the $SU(3)$ orbits of integral curves of $E_4$.

\end{itemize}
\section*{Appendix A}
\appendix
\setcounter{equation}{0}
\def\theequation{\thesection{A}\arabic{equation}}
In this Appendix we establish a non--existence result which narrows down the classification of integrable backgrounds for conformal geodesic equations.
\begin{prop}
There are no non-flat Riemannian Gibbons-Hawking metrics with three commuting Killing vectors.
\end{prop}
\noindent
{\bf  Proof.}
Suppose there are, and we have the Gibbons-Hawking metric (\ref{GH}) with two further Killing vectors, say $L_1, L_2$, which commute with each other and with $K=\partial_\tau$. We find restrictions on the $L_1, L_2$ that force flatness. First note that for either $L$
\[L^a\partial_a(V^{-1})=L^a\partial_a(K^bK_b)=2L^aK^b\nabla_aK_b=2K^aK^b\nabla_aL_b,\]
using $[K,L]=0$, and this is manifestly zero, so that $V$ is constant along each $L$. 
Decompose each $L$ as
\[L^a\partial_a=L^0\partial_\tau+L^i\partial_i,\]
then the vanishing of the commutator with $K$ gives
\[\partial_\tau L^0=0=\partial_\tau L^i.\] Next with
\[h_{ab}:=g_{ab}-VK_aK_b=V\delta_{ij},\]
where $\delta_{ij}=0$ unless $i,j=1,2,3$ when it is the usual Kronecker delta, we have
\[\mathcal{L}_Lh_{ab}=0,\]
so that $L^i\partial_i$ is a Killing vector of the flat (Riemannian) 3-metric $\delta_{ij}$. Next, the $(i,j)$ component of
\[0=\mathcal{L}_Lg_{ab}=L^a\partial_ag_{bc}+g_{ac}\partial_bL^a+g_{ba}\partial_cL^a\]
forces $L^0$ to be constant when it can be set to zero by subtracting a constant multiple of $K$ from $L$. Finally from the $b=j$ component of
\[0=\mathcal{L}_LK_b\]
we see that 
\[0=L^i\partial_i(V^{-1}\omega_j)+V^{-1}\omega_j\partial_iL^j,\]
so that $L^i$ Lie-drags $\omega$.

Thus any further Killing vector is w.l.o.g. a symmetry of Riemannian flat 3-space, Lie-dragging both $V$ and $\omega$. These can be translations or rotations. Two independent rotations will not commute so $L_1, L_2$ can be either two independent translations, in which case we may suppose them orthogonal, or a translation and a rotation, in which case commutation requires them to be orthogonal. There are two cases to consider:
\begin{itemize}
 \item Suppose $L_1=\partial_x, L_2=\partial_y$ then up to constants $V=z$ or $V=1$. Now if $V\neq 1$
 \[d\omega=*dV=dx\wedge dy,\]
 and there is no $\omega$ preserved by both $L_1$ and $L_2$. Thus $V=1$ when $\omega=0$ and this class of solutions is flat.
 \item Suppose $L_1=\partial_z$ and $L_2=\partial_\phi$ in cylindrical polars. Then $V$ depends only on $\rho$ and is harmonic so up to constants $V=\log\rho$ or $V=1$. Now if $V\neq 1$
 \[d\omega=\frac{*d\rho}{\rho}=dz\wedge d\phi,\]
 and again there is no acceptable 
solution\footnote{
Note that the claim is not true in the neutral signature:
Consider the simplest ASD pp--wave
\[
g=dwdx+dzdy+y^2dw^2,
\]
and complete the  square to put it in the Gibbons--Hawking form
\[
g=V\Big(dYdz-\frac{1}{4}dx^2\Big)+V^{-1}\Big(dw+\frac{1}{2}(3Y)^{-2/3}dx
\Big)^2,
\]
where $Y=y^3/3$ and $V=(3Y)^{-2/3}$.
This metric has a 9D isometry group \cite{CDT} which is nilpotent, and
includes the Abelian 3D algebra
spanned by $(\p_x, \p_w, \p_z)$.}, except constant $V$.
\end{itemize}
\koniec
\section*{Appendix B}
\appendix
\setcounter{equation}{0}
\def\theequation{\thesection{B}\arabic{equation}}
In this appendix we shall give necessary and sufficient conditions
for a trajectory of a Killing vector in a four--dimensional Einstein space to be a conformal
geodesic
\begin{prop}
Let $K$ be a non--null Killing vector on a (pseudo) Riemannian Einstein 
four manifold $(M, g)$. The trajectories of $K$ are conformal geodesics
if and only if
\be
\label{paul_criteron}
d |K|^2\wedge  W=0, \quad \mbox{where}\quad W=\star_g(K\wedge dK).
\ee
\end{prop}
\noindent
{\bf Proof.}
Suppose the Killing vector is $K$ with norm $V^2=g(K,K)$. We claim that, for any Killing vector
\[\nabla_aK_b=\frac{2}{V}V_{[a}K_{b]}+\frac{1}{2V^2}\epsilon_{ab}^{\;\;\;\;\;\;cd}W_cK_d,\]
for some $W_a$ orthogonal to $K$; as a form $W=*(K\wedge dK)$. Now the velocity will be $u^a=V^{-1}K^a$ and the acceleration as above will be $a_b=-V^{-1}\nabla_bV$ so we need
\[-|a|^2u_b=u^c\nabla_c{a}_b=V^{-1}K^c\nabla_c(-V^{-1}V_b)=-V^{-2}K^c\nabla_bV_c=V^{-2}V^c\nabla_bK_c=
\]
\[
=V^{-2}V^c(\frac{2}{V}V_{[b}K_{c]}+\frac{1}{2V^2}\epsilon_{bc}^{\;\;\;\;\;\;de}W_dK_e),
\]
i. e.
\[-V^{-3}|\nabla V|^2K_b=-V^{-3}|\nabla V|^2K_b+\frac12V^{-4}\epsilon_b^{\;\;cde}V_cW_dK_e,\]
so the last term needs to vanish, which is equivalent to
(\ref{paul_criteron}).
\koniec
As a special case we recover the result of \cite{tod_circles}, where it was 
shown that if $W=0$, then the Killing trajectories are conformal geodesics.
We also find that in the case of the Fubini--Study metric 
(\ref{cp2metric}) on $\CP^2$ the trajectories of
$\p/\p \psi$ are conformal geodesics, but the trajectories of
$K=\p/\p\phi$ (or any other Killing vector) are not. 
\section*{Appendix C}
\appendix
\setcounter{equation}{0}
\def\theequation{\thesection{C}\arabic{equation}}
In this appendix we shall reformulate the Lorentz force equations, as well as
some of the calculations underlying the conformal Killing--Yano tensors in terms of the two component
spinor notation \cite{PR}. 

The Hodge $\ast$ operator of an oriented Riemannian manifold $(M, g)$ is an
involution on two-forms, and induces a decomposition
\[
\Lambda^{2}(T^*M) = \Lambda_{+}^{2}(T^*M) \oplus \Lambda_{-}^{2}(T^*M)
\]
of two-forms into self-dual (SD)
and anti-self-dual (ASD)  components. 
Locally there exist complex rank-two vector bundles $\spp, \spp'$  over $M$ equipped with covariantly constant symplectic structures
$\epsilon, \epsilon'$ such that
\[
\C\otimes TM\cong \spp\otimes \spp', \quad  g=\epsilon\otimes\epsilon',
\]
and
\[
\Lambda_{+}^{2}\cong {\spp'}^*\odot {\spp'}^*, \quad
\Lambda_{-}^{2}\cong {\spp}^*\odot {\spp}^*.
\]
Elements of $\spp$ are of the form $\kappa^A=(\kappa^0, \kappa^1)$, and the spinor indices $A, B, C, \dots =0, 1$ are lowered
using the symplectic form $\epsilon_{AB}$ with $\epsilon_{01}=1$. 
Therefore any two-form admits a spinor decomposition
\[
F_{ab}=f_{AB}\epsilon_{A'B'}+\tilde{f}_{A'B'}\epsilon_{AB},
\]
where $f_{AB}=f_{(AB)}$, and $\tilde{f}_{A'B'}=\tilde{f}_{(A'B')}$. In particular the parallel
basis $(\Omega^1, \Omega^2, \Omega^3)$ of $\Lambda_{+}^2$ on a hyper--K\"ahler
four manifold decomposes as
\[
{\Omega^i}_{ab}={\phi^i}_{A'B'}\epsilon_{AB}, \quad i=1, 2, 3.
\]
\vskip5pt
Consider the Lorentz force equation (\ref{lorentzf}) with self--dual Maxwell field $F$, and set
\[
F_{ab}=|a|\phi_{A'B'}\epsilon_{AB}
\]
where $\phi_{A'B'}=\phi_{(A'B')}$ is a valence two symmetric spinor.
The Lorentz force equation takes the form
\[
u^{CC'}\nabla_{CC'} u_{AA'}=-|a| \phi_{A'B'} {u^{B'}}_A
\]
and $\phi_{A'B'}\phi^{A'B'}=2$. The ASD two--form $W$ from (\ref{Wtensor}) is
\[
W_{ab}=\omega_{AB}\epsilon_{A'B'},
\]
and the CKY equation (\ref{CKY}) takes the form
\be
\label{ks1}\nabla_{AA'}\omega_{BC}=\epsilon_{A(B}K_{C)A'},\quad\mbox{where}\quad K=\p/\p\psi.
\ee
First note that $\omega_{AB}$ gives rise to three first integrals of the geodesic
equations. To see this consider the quadratic expression in velocity $\omega_{AB}{\phi^i}_{A'B'}u^{AA'}u^{BB'}$, whose 
rate of change is
\[
u^{CC'}\nabla_{CC'}(\omega_{AB}{\phi^i}_{A'B'}u^{AA'}u^{BB'})=\frac12u^{AA'}
{\phi^i}_{A'B'}K_{A}^{\;\;B'}=\frac12u^{AA'}\nabla_{AA'}x^i,\]
by (\ref{ks1}) and (\ref{h1}) so that three constants of the motion for the geodesic equation are defined by
\[q^i:=\omega_{AB}{\phi^i}_{A'B'}u^{AA'}u^{BB'}-\frac12x^i.
\]

It is known from \cite{tod_circles} that a Killing spinor will give constants of the motion for the conformal geodesic equation: consider $\omega_{A'B'}u^{AA'}a_A^{\;\;B'}$ then its rate of change following (
\ref{conf_circ}) and the definition $a=\nabla_u u$ is
\[u^{CC'}\nabla_{CC'}(\omega_{AB}u^{AA'}a_{A'}^{\;\;B})=-\frac12K_ba^b=-
\frac12 u^c\nabla_c(K_bu^b),\]
so that the following is a conserved quantity for the conformal geodesic equation:
\[\omega_{AB}u^{AA'}a_{A'}^{\;\;B}+\frac12K_bu^b.\]
In the present case, when we have the relation
\[
u_{AA'}{a^{A}}_{B'}=u_{A(A'}{a^{A}}_{B')}=-\frac{1}{2}|a|\phi_{A'B'}
\]
for a real constant spinor $\phi_{A'B'}$
between velocity and acceleration, this constant becomes
\be\label{ks2}-|a|u^{AA'}u^{BB'}\phi_{A'B'}\omega_{AB}+\frac12K_bu^b.\ee
This is the first integral ${\mathcal W}$ in (\ref{first_ints}).

\end{document}